\documentclass[a4paper,11pt]{article}
\usepackage[margin=28mm]{geometry}
\usepackage[T1]{fontenc}
\usepackage[utf8]{inputenc}
\usepackage{amsmath,amssymb,amsfonts,amsthm}
\usepackage{graphicx,xcolor,tikz}
\usetikzlibrary{arrows,arrows.meta,positioning,calc}
\usepackage[font=small,labelfont=bf]{caption}
\usepackage[hidelinks]{hyperref}
\newcommand{\invlim}{\varprojlim}
\newcommand{\id}{\operatorname{id}}

\definecolor{fiberblue}{RGB}{35,100,165}
\definecolor{branchgreen}{RGB}{0,125,110}
\tikzset{maparrow/.style={-{Stealth[length=2mm]},semithick},
axis/.style={-{Stealth[length=1.5mm]},thin,gray},
fiber/.style={fiberblue,line width=2.5pt},
graphline/.style={fiberblue,thick}}

\newtheorem{theorem}{Theorem}[section]

\newtheorem{definition}[theorem]{Definition}
\newtheorem{lemma}[theorem]{Lemma}

\newtheorem{example}[theorem]{Example}

\newtheorem{corollary}[theorem]{Corollary}

\newtheorem{observation}[theorem]{Observation}

\newcommand{\diam} {\mathop{\rm diam}\nolimits}

\newcommand{\Cl}  {\mathop{\rm Cl}\nolimits}
\newcommand{\Int}  {\mathop{\rm Int}\nolimits}

\begin{document}

\title{Universal and Reciprocating Domains in Classes of Continua and Beyond}
\date{}
\author{Iztok Bani\v c, Goran Erceg, Judy Kennedy, Chris Mouron, Van Nall}
\maketitle

\begin{abstract}
We introduce and study the notions of reciprocating domains in classes of topological spaces. A reciprocating domain is a member of a class that maps onto every member which maps onto it. We provide general characterizations, showing in particular that in upward directed classes reciprocating domains coincide with universal domains. We investigate these domains in several classes of compact metrizable spaces and continua. We also construct a chainable universal domain admitting a continuous monotone surjection onto a pseudoarc although every nondegenerate subcontinuum of the domain contains an arc. In particular, a chainable reciprocating domain need not contain a pseudoarc. 
\end{abstract}
\-

\section{Introduction}

A universal domain for a class of topological spaces is a member of
that class which maps continuously onto every other member. Classical
examples include the Cantor set for non-empty compact metrizable spaces
and the arc for Peano continua. A reciprocating domain satisfies a
weaker condition: whenever a member of the class maps continuously
onto it, a continuous surjection also exists in the other direction.
These notions describe, respectively, maximum and maximal elements
of the ordering induced by continuous surjections, after spaces that
map onto one another are identified.

Continuous surjections and embeddings give different ways of comparing
continua. The pseudoarc illustrates the distinction: it maps onto every
chainable continuum, whereas each of its nondegenerate subcontinua is
again a pseudoarc. We show that a chainable continuum can map onto a
pseudoarc and still contain no pseudoarc. More precisely, we construct
a chainable universal domain in which every nondegenerate subcontinuum
contains an arc. This construction complements the general theory of
reciprocating domains by showing that a maximal mapping property need
not determine the hereditary structure of the domain.

Throughout the paper, all topological spaces are assumed to be
non-empty. A continuum is a non-empty compact connected metrizable
space. We allow a continuum to consist of one point. We first introduce
the two types of domains and their order-theoretic interpretations.

\begin{definition}
	Let $\mathcal X$ be a class of topological spaces and let $X\in \mathcal X$. If for each $Y\in \mathcal X$, there is a continuous surjection from $X$ to $Y$, then we say that $X$ is a \emph{universal domain in} $\mathcal X$. 
\end{definition}

\begin{observation}\label{nud}
Let $\mathcal T$ be the class of all topological spaces. Note that for any $X\in \mathcal T$ there is a $Y\in \mathcal T$ such that the cardinality of $Y$ is greater than the cardinality of $X$. It follows that for any $X\in \mathcal T$ there is a $Y\in \mathcal T$ such that there are no continuous surjections from $X$ to $Y$. Therefore, there are no universal domains in $\mathcal T$.
\end{observation}
\begin{observation}\label{MCcantor}
		Let $\mathcal X$ be the class of compact metrizable spaces and let $C$ be a Cantor set. 
Then $C\in \mathcal X$. Since every non-empty compact metrizable space is a continuous image of $C$, 
	it follows that the Cantor set $C$ is a universal domain in $\mathcal{X}$.
\end{observation}
\begin{observation}
	Let $\mathcal X$ be the class of Peano continua and let $I=[0,1]$. 	Then $I\in \mathcal X$.
	Since every Peano continuum is a continuous image of $I$, it follows that $I$ is a universal domain in $\mathcal{X}$.
\end{observation}
In the following observation, we present an interpretation of the notion of  universal domains. 
\begin{observation}\label{milan}
	Let $\mathcal X$ be a class of topological spaces. Then one can define an equivalence relation $\sim$ on the class $\mathcal{X}$ by
	$$
	X\sim Y ~~~ \Longleftrightarrow ~~~ \textup{there are  continuous surjections } f:Y\rightarrow X \textup{ and } g:X\rightarrow Y
 	$$ 
	for all $X,Y\in \mathcal{X}$. For each $X\in \mathcal{X}$, we also define its equivalence class as
	$$
	[X]=\{Y\in \mathcal{X} \mid Y\sim X\}.
	$$
We denote by $\mathcal{X}/_{\sim}$ the collection
	$$
	\mathcal{X}/_{\sim}= \{[X] \mid X\in \mathcal{X}\}
	$$
	of classes of topological spaces. On this collection, we define the ordering relation $\leq$ as follows: for all $X,Y\in \mathcal{X}$,
	$$
	[X] \leq [Y] ~~~ \Longleftrightarrow ~~~ \textup{there is a  continuous surjection } f:Y\rightarrow X.
	$$
	Then the relation $\leq$ is an order on $\mathcal{X}/_{\sim}$ and the pair
	$$
	(\mathcal{X}/_{\sim}, \leq)
	$$
	may be considered as a partially ordered collection of classes of topological spaces. In this partially ordered collection, we interpret the meaning of  equivalence classes of universal domains. Let $X\in \mathcal X$. Then the following statements are equivalent:
		\begin{itemize}
			\item $X$ is a universal domain in $\mathcal X$.
			\item For each $Y\in \mathcal X$, there is a continuous surjection from $X$ to $Y$.
			\item For each $Y\in \mathcal X$, $[Y]\leq [X]$.
			\item $[X]$ is the maximum element in $(\mathcal{X}/_{\sim}, \leq)$.
		\end{itemize}
		This proves that $X$ is a universal domain in $\mathcal X$ if and only if $[X]$ is the maximum element in $(\mathcal{X}/_{\sim}, \leq)$.  
		\end{observation}
		So, maxima in $(\mathcal{X}/_{\sim}, \leq)$ are exactly the classes of topological spaces that consist of universal domains in $\mathcal X$.  In light of this order-theoretic interpretation of universal domains, it is only natural to investigate also the maximal elements in $(\mathcal{X}/_{\sim},\leq)$ and to study their representatives within the class $\mathcal X$. The following observation provides such a natural characterization of maximal elements in $(\mathcal{X}/_{\sim},\leq)$.	
		\begin{observation}
			Let $\mathcal X$ be a class of topological spaces and let $X\in \mathcal X$. Then the following statements are equivalent:
		\begin{itemize}
			\item $[X]$ is a maximal element in $(\mathcal{X}/_{\sim}, \leq)$.
			\item For each $Y\in \mathcal X$, it holds that if $[X]\leq [Y]$, then $[Y] = [X]$.
			\item For each $Y\in \mathcal X$, it holds that if $[X]\leq [Y]$, then $[Y]\leq [X]$.
			\item For each $Y\in \mathcal X$, it holds that if there is a continuous surjection from $Y$ to $X$, then there is a continuous surjection from $X$ to $Y$.
		\end{itemize}
		\end{observation}	
		Inspired by the above observation, we introduce a reciprocating domain of a class of topological spaces.
\begin{definition}
	Let $\mathcal X$ be a class of topological spaces and let $X\in \mathcal X$. We say that $X$ is a \emph{reciprocating domain in $\mathcal X$} if for each $Y\in\mathcal X$, the following holds:
	$$
	\textup{there is a continuous surjection } Y \to X \Longrightarrow   \textup{ there is a continuous surjection } X \to Y.
	$$
\end{definition}
\begin{observation}\label{karin}
	Let $\mathcal X$ be a class of topological spaces and let $X\in \mathcal X$. Note that if $X$ is a universal domain in $\mathcal X$, then $X$ is a reciprocating domain in $\mathcal X$.
\end{observation}
\begin{example}
		Let $\mathcal X$ be the class of compact metrizable spaces and let $C$ be a Cantor set. 
Then $C$ is a reciprocating domain in $\mathcal{X}$.
\end{example}
\begin{example}
	Let $\mathcal X$ be the class of Peano continua and let $I=[0,1]$. 	Then $I$ is a reciprocating domain in $\mathcal{X}$.
\end{example}
\begin{example}\label{psajd}
		Let $X,Y\in \mathcal X$ be such that there are no continuous surjections from $X$ to $Y$ and there are no continuous surjections from $Y$ to $X$. Such spaces $X$ and $Y$ do exist; for example, among other things,  Waraszkiewicz showed in \cite{W} that there are such spaces $X$ and $Y$. Let $\mathcal X=\{X,Y\}$. Then both, $X$ and $Y$ are reciprocating domains in $\mathcal{X}$ but there is no universal domain in $\mathcal X$.
\end{example} 

Recent nonexistence results also show how restrictions on the
members of a class affect the existence of a common domain. In
\cite{Matsuhashi2026}, Matsuhashi studies \emph{common models},
that is, continua that map continuously onto every member of a given
class, without requiring the model to belong to that class. Using
meandering continua, he proves that, for each continuum $Z$, the
class of hereditarily indecomposable continua admitting a continuous
light surjection onto $Z$ has no common model. He also obtains
nonexistence results for classes arising from failures of Whitney
reversibility. These results concern the same mapping relation as
universal domains and show that nonexistence can persist even when
the common domain is allowed to lie outside the class.

Alongside its new results, this paper provides an overview of
universal domains in several natural classes of continua and compact
metrizable spaces. The general characterization in
Theorem~\ref{cilka}, the smooth-fan characterization in
Theorem~\ref{SF}, and the chainable domain containing no pseudoarc
in Theorem~\ref{pa-thm:counterexample} are central parts of the
paper. The common indecomposable-domain construction in
Theorem~\ref{ired} also gives a direct proof that the class of
indecomposable continua is directed upwards. These results are
presented together with classical universality theorems: the Cantor
set for compact metrizable spaces, the Hahn--Mazurkiewicz theorem
for Peano continua, the pseudoarc for weakly chainable continua,
and the Cantor fan for uniformly path connected continua and
smooth dendroids. The nonexistence results of Waraszkiewicz,
Ingram, and Russo provide the corresponding background for classes
without universal domains. In this way, the examples both develop
the theory of reciprocating domains and collect the existing
mapping results on which that theory rests.

The universality of the pseudoarc among chainable continua was
proved by Mioduszewski in \cite{Mioduszewski}; another proof was
given by Irwin and Solecki in \cite{IrwinSolecki}. Continuous images
of the pseudoarc were studied by Lelek and Fearnley in
\cite{Lelek,Fearnley}. We use these classical mapping results as
background for the construction below.

The chainable case also leads to a question about subcontinua. A
surjection onto a pseudoarc forces the domain to contain, over every
composant of the range, a strictly decreasing sequence of nondegenerate
indecomposable subcontinua; see Theorem~\ref{bujaha}. Since distinct
composants of a pseudoarc are disjoint and there are uncountably many
of them, the domain contains uncountably many pairwise disjoint
nondegenerate indecomposable subcontinua. Nevertheless, the domain
need not contain a pseudoarc. We prove this in
Theorem~\ref{pa-thm:counterexample} by replacing points of compatible
countable dense sets in the coordinate intervals of an inverse sequence
by nondegenerate intervals and lifting its bonding maps. The induced
map is a continuous monotone surjection, and every nondegenerate
subcontinuum of its domain contains an arc.

The replacement of points by continua is related to condensation of
singularities, developed for arc-like continua by Ma\'{c}kowiak in
\cite{Mackowiak}. Krzempek used this technique to prescribe spaces of
end points of chainable continua in \cite{Krzempek}; his Theorem 4
recalls an atomic replacement theorem of Pol and Re\'{n}ska \cite[Theorem 3.2]{PolRenska} for
zero-dimensional closed sets. Here, the replacements are made in the
coordinate intervals, and chainability follows from the resulting
inverse sequence. The arc construction also has precedents: Toledo
used homeomorphic restrictions of bonding maps to obtain an arc in
\cite[Theorem 4.12]{Toledo}, while Charatonik, Mena, and Roe proved in
\cite[Proposition 4.13]{CharatonikMenaRoe} that every nondegenerate
subcontinuum of each of six specified inverse-limit spaces contains
an arc. We apply a related argument within the fibers of the induced
map, and then extract a nondegenerate subcontinuum in a fiber from an
arbitrary nondegenerate subcontinuum of the domain.

\section{When reciprocating domains are universal}\label{sec:general}
In this section, we study the relationship between reciprocating and universal domains. We characterize classes in which these notions coincide and apply our results to disjoint unions, wedges, and the class of all topological spaces.
\begin{theorem}\label{general:1}
	Let $\mathcal X$ be a class of topological spaces,  let $X,Y\in \mathcal X$ and let $f:Y\to X$ be a continuous surjection. If $X$ is a reciprocating domain in $\mathcal X$, then also $Y$ is a reciprocating domain in $\mathcal X$.
\end{theorem}
\begin{proof}
	Suppose that $X$ is a reciprocating domain in $\mathcal X$ and let $Z\in\mathcal X$ be such that there is a continuous surjection from $Z$ to $Y$. Let $g:Z\to Y$ be such a surjection. Since $X$ is a reciprocating domain in $\mathcal X$ and $f\circ g:Z\to X$ is surjective, there is a continuous surjection from $X$ to $Z$. Let $h:X\to Z$ be such a surjection. Then $h\circ f:Y\to Z$ is a continuous surjection. Therefore,  $Y$ is a reciprocating domain in $\mathcal X$.
\end{proof}

\begin{theorem}\label{zoesita}
	Let $\mathcal X$ be a class of topological spaces and let $X\in \mathcal X$. The following statements are equivalent.
	\begin{enumerate}
		\item\label{aa1} $X$ is a universal domain in $\mathcal{X}$.
		\item\label{aa2} $X$ is a reciprocating domain in $\mathcal{X}$ and there is a universal domain in $\mathcal X$.
	\end{enumerate}
\end{theorem}

\begin{proof}
	To prove the implication from \ref{aa1} to \ref{aa2}, suppose that $X$ is a universal domain in $\mathcal X$. Let $Y\in \mathcal{X}$ be any space such that there is a continuous surjection from $Y$ to $X$. Since $X$ is a universal domain in $\mathcal X$, there is a continuous surjection from $X$ to $Y$, as required. Also, since $X$ is a universal domain in $\mathcal{X}$, it follows that there is a universal domain in $\mathcal X$.

	To prove the implication from \ref{aa2} to \ref{aa1}, suppose that $X$ is a reciprocating domain in $\mathcal{X}$ and let $Y$ be a universal domain in $\mathcal X$. To prove that $X$ is a universal domain in $\mathcal{X}$, let $Z\in \mathcal{X}$. We prove that there is a continuous surjection from $X$ to $Z$. Since $Y$ is a universal domain in $\mathcal X$, there are continuous surjections $h:Y\to Z$ and $g:Y\rightarrow X$. By the assumption that $X$ is a reciprocating domain in $\mathcal{X}$, this implies that there is also a continuous surjection $f$ from $X$ to $Y$. Therefore, $h\circ f$ is a continuous surjection from $X$ to $Z$. Hence, statement \ref{aa1} holds.
\end{proof} 

\begin{definition}
	Let $\mathcal X$ be a class of topological spaces. If for all $X,Y\in \mathcal X$ there is $Z\in \mathcal X$ such that there are continuous surjections from $Z$ to $X$ and from $Z$ to $Y$, then we say that the class $\mathcal X$ is \emph{directed upwards}. 
\end{definition}
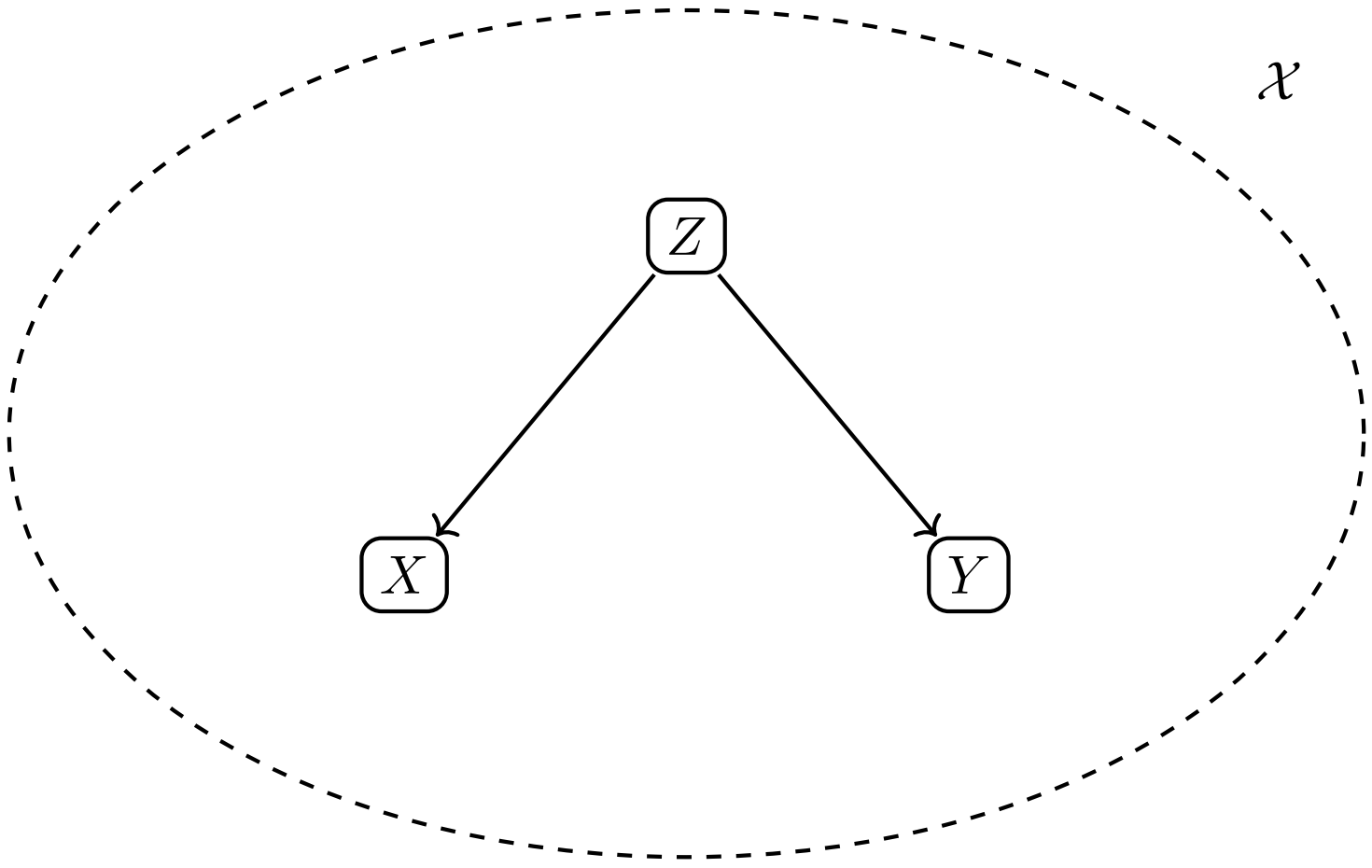
\begin{center}
\begin{tikzpicture}[->,thick]

\draw[dashed] (0,0) ellipse (4.8 and 3);
\node at (4.2,2.5) {$\mathcal X$};

\node[draw,rounded corners] (Z) at (0,1.4) {$Z$};
\node[draw,rounded corners] (X) at (-2,-1) {$X$};
\node[draw,rounded corners] (Y) at (2,-1) {$Y$};

\draw (Z) -- (X);
\draw (Z) -- (Y);

\end{tikzpicture}
\end{center}
\begin{theorem}\label{cilka}
	Let $\mathcal X$ be a class of topological spaces which is directed upwards and let $X\in \mathcal X$. Then
	$$
	X \textup{ is a reciprocating domain in } \mathcal{X} ~~~ \Longleftrightarrow ~~~ X \textup{ is a universal domain in } \mathcal{X}.
	$$
\end{theorem}
\begin{proof}
	Suppose that $X$ is a reciprocating domain in $\mathcal{X}$ and that it is not a universal domain in $\mathcal X$. Let $Y\in \mathcal X$ be such that there are no continuous surjections from $X$ to $Y$. Also, let $Z\in \mathcal X$ be such that there are continuous surjections from $Z$ to $X$ and from $Z$ to $Y$. Since $X$ is a reciprocating domain in $\mathcal X$, it follows that there is a continuous surjection from $X$ to $Z$. Let $f:X\to Z$ be a continuous surjection and let $g:Z\to Y$ be a continuous surjection. Then $g\circ f$ is a continuous surjection from $X$ to $Y$, which is a contradiction. Therefore,   $X$ is a universal domain in $\mathcal{X}$. The other implication follows from Observation \ref{karin}. 
\end{proof}
\begin{theorem}
	Let $\mathcal X$ be a class of topological spaces such that there is a reciprocating domain in $\mathcal X$. The following statements are equivalent.
	\begin{enumerate}
		\item\label{muco11} For each $X\in \mathcal X$,
	$$
	X \textup{ is a reciprocating domain in } \mathcal{X} ~~~ \Longleftrightarrow ~~~ X \textup{ is a universal domain in } \mathcal{X}.
	$$
	\item\label{muco22} The class $\mathcal X$ is directed upwards.
	\end{enumerate}
\end{theorem}
\begin{proof}
	Let $X_0$ be a reciprocating domain in $\mathcal X$. First, we prove the implication from \ref{muco11} to \ref{muco22}. Suppose that for each $X\in \mathcal X$,
	$$
	X \textup{ is a reciprocating domain in } \mathcal{X} ~~~ \Longrightarrow ~~~ X \textup{ is a universal domain in } \mathcal{X}.
	$$
	It follows that $X_0$ is a universal domain. To prove that $\mathcal X$ is directed upwards, let $X,Y\in \mathcal X$ and let $Z=X_0$. Then (since $X_0$ is a universal domain) there are continuous surjections from $Z$ to $X$ and from $Z$ to $Y$. Therefore, $\mathcal X$ is directed upwards. The implication from \ref{muco22} to \ref{muco11} follows from Theorem \ref{cilka}. 
\end{proof}
In Theorem \ref{unijeunije22}, we demonstrate how Theorem \ref{cilka} may be applied to show that there are no reciprocating domains in the class $\mathcal T$ of all topological spaces. We use simple results about disjoint unions (Theorem \ref{unijeunije1}) and wedges (Theorem \ref{wedge}) of topological spaces to do so. We start with  disjoint unions of topological spaces
\begin{definition}
	Let $X$ and $Y$ be topological spaces, let $\mathcal T_X$ be the topology on $X$ and let $\mathcal T_Y$ be the topology on $Y$. Then we define the topology on the disjoint union $X\sqcup Y$ to be the topology with the basis $\mathcal T_X\cup \mathcal T_Y$. 
\end{definition}
\begin{theorem}\label{unijeunije}
	Let $\mathcal X$ be a class of topological spaces such that for all $X,Y\in \mathcal X$, there is $Z\in \mathcal X$ such that $Z$ is homeomorphic to $X\sqcup Y$. Then $\mathcal X$ is directed upwards. 
\end{theorem}
\begin{proof}
	Let $X,Y\in \mathcal X$, let $Z\in \mathcal X$ be such that $Z$ is homeomorphic to $X\sqcup Y$, let $x_0\in X$ and let $y_0\in Y$. Also, let $\varphi:Z\rightarrow X\sqcup Y$ be a homeomorphism. Note that $f:Z\rightarrow X$ and $g:Z\to Y$, defined by  
	$$
	f(x)=\begin{cases}
				\varphi(x)\text{;} & \varphi(x)\in X\\
				x_0\text{;} & \varphi(x)\in Y
			\end{cases} ~~~ \textup{ and } ~~~ g(x)=\begin{cases}
				\varphi(x)\text{;} & \varphi(x)\in Y\\
				y_0\text{;} & \varphi(x)\in X
			\end{cases}
	$$
	for each $x\in Z$, are continuous surjections from $Z$ to $X$ and from $Z$ to $Y$. Therefore, $\mathcal X$ is directed upwards.  
\end{proof}
\begin{corollary}\label{unijeunije1}
	Let $\mathcal X$ be a class of topological spaces such that for all $X,Y\in \mathcal X$, there is $Z\in \mathcal X$ such that $Z$ is homeomorphic to $X\sqcup Y$, and let $X_0\in \mathcal X$. If $X_0$ is a reciprocating  domain in $\mathcal{X}$, then $X_0$ is a universal domain in $\mathcal{X}$. 
\end{corollary}
\begin{proof}
	Suppose that $X_0$ is a reciprocating  domain in $\mathcal{X}$. Since $\mathcal X$ is directed upwards by Theorem \ref{unijeunije}, it follows from Theorem \ref{cilka} that $X_0$ is a universal domain in $\mathcal{X}$. 
\end{proof}
We continue with  wedges of topological spaces. They are defined with topological quotients.

\begin{definition}
	Let $X$ be a topological space and let $\sim$ be an equivalence relation on $X$. For each $x\in X$, we use $[x]$ to denote the equivalence class of the element $x$ with respect to the relation $\sim$. 
	We also use $X/_{\sim}$ to denote the quotient space $X/_{\sim}=\{[x] \ | \ x\in X\}$. 
	Also,  let $q:X\rightarrow X/_{\sim}$ be the quotient map that is defined by $q(x)=[x]$ for each $x\in X$,  and let $U\subseteq X/_{\sim}$. Then 
	$$
	U \textup{ is open in } X/_{\sim} ~~~ \Longleftrightarrow ~~~ q^{-1}(U)\textup{ is open in } X.
	$$ 
\end{definition}
\begin{definition}
	Let $X$ be a topological space, let $A\subseteq X$ and let $\sim$ be the equivalence relation on $X$ that is defined as follows:
	$$
	x\sim y ~~~ \Longleftrightarrow ~~~ x=y ~~~ \textup{ or } ~~~ x,y\in A.
	$$
	Then we use $X/_A$ to denote the quotient space $X/_{\sim}$.
\end{definition}
\begin{definition}
	Let $X$ and $Y$ be topological spaces, let $x\in X$ and let $y\in Y$. We say that the quotient space $(X\sqcup Y)/_{\{x,y\}}$ is a \emph{wedge of $X$ and $Y$}. 
\end{definition}
\begin{theorem}\label{wedge}
	Let $\mathcal X$ be a class of topological spaces such that for all $X,Y\in \mathcal X$, there is $Z\in \mathcal X$ such that for some wedge $W$ of $X$ and $Y$, $Z$ is homeomorphic to $W$. Then $\mathcal X$ is directed upwards.  
\end{theorem}
\begin{proof}
Let $X,Y\in\mathcal X$, choose $x_0\in X$ and $y_0\in Y$,
and let $W=(X\sqcup Y)/_{\{x_0,y_0\}}$ be a wedge homeomorphic
to a member $Z$ of $\mathcal X$. The maps on $X\sqcup Y$ which
are the identity on $X$ and constant with value $x_0$ on $Y$,
and the identity on $Y$ and constant with value $y_0$ on $X$,
respectively, are continuous and constant on the identified pair.
They therefore induce continuous surjections $W\to X$ and
$W\to Y$. Composing with a homeomorphism $Z\to W$ proves that
$\mathcal X$ is directed upwards.
\end{proof}
\begin{corollary}\label{wedge1}
	Let $\mathcal X$ be a class of topological spaces such that for all $X,Y\in \mathcal X$, there is $Z\in \mathcal X$ such that for some wedge $W$ of $X$ and $Y$, $Z$ is homeomorphic to $W$, and let $X_0\in \mathcal X$. If $X_0$ is a reciprocating  domain in $\mathcal{X}$, then $X_0$ is a universal domain in $\mathcal{X}$.
\end{corollary}
\begin{proof}
	Suppose that $X_0$ is a reciprocating  domain in $\mathcal{X}$. Since $\mathcal X$ is directed upwards by Theorem \ref{wedge}, it follows from Theorem \ref{cilka} that $X_0$ is a universal domain in $\mathcal{X}$. 
\end{proof}
We use Corollary \ref{unijeunije1} (or Corollary \ref{wedge1}) to show that there are no reciprocating domains in the class of all topological spaces $\mathcal T$.
\begin{theorem}\label{unijeunije22}
	Let $\mathcal T$ be the class of all topological spaces. For each $X\in \mathcal T$, $X$ is not a reciprocating domain in $\mathcal T$.
\end{theorem}
\begin{proof}
	Suppose that there is $X_0\in \mathcal T$ such that $X_0$ is a reciprocating domain in $\mathcal T$. Since for all $X,Y\in \mathcal T$, there is $Z\in \mathcal X$ such that $Z$ is homeomorphic to $X\sqcup Y$ (or such that for some wedge $W$ of $X$ and $Y$, $Z$ is homeomorphic to $W$), it follows from Corollary \ref{unijeunije1} (or Corollary \ref{wedge1}) that $X_0$ is a universal domain in $\mathcal T$, which is a contradiction with Observation \ref{nud}. 
\end{proof}
We conclude the section with the following theorem. 

\begin{theorem}\label{zelosplosno}
	Let $\mathcal X$ be a class of topological spaces, let $U\in \mathcal X$ be a universal domain in $\mathcal X$ and let $X\in \mathcal X$. The following statements are equivalent.
	\begin{enumerate}
		\item\label{juc1} $X$ is a reciprocating domain in $\mathcal X$.
		\item\label{juc2} There is a continuous surjection $f:X\to U$.
		\item\label{juc3} $X$ is a universal domain in $\mathcal X$.
	\end{enumerate}  
\end{theorem}
\begin{proof}
	First, we prove the implication from \ref{juc1} to \ref{juc2}. Suppose that $X$  is a reciprocating domain in $\mathcal X$. Since $U$ is a universal domain in $\mathcal X$, it follows that there is a continuous surjection from $U$ to $X$. Since $X$ is a reciprocating domain in $\mathcal X$, there is also a continuous surjection from $X$ to $U$.  
	To prove the implication from \ref{juc2} to \ref{juc1}, suppose that there is a continuous surjection $f:X\to U$, let $Y\in \mathcal X$ and suppose that there is a continuous surjection from $Y$ onto $X$.  Since $U$ is a universal domain in $\mathcal X$, it follows that there is a continuous surjection $g:U\to Y$. Let $g$ be such a function. Then $g\circ f$ is a continuous surjection from $X$ to $Y$. Therefore,  $X$ is a reciprocating domain in $\mathcal X$. 
	The equivalence of \ref{juc1} and \ref{juc3} follows from Theorem \ref{zoesita}. 
\end{proof}

\section{Examples in certain classes of metrizable spaces}
In this section, we apply the general results of the previous section to several classes of metrizable spaces. We characterize reciprocating domains among compact metrizable spaces, Peano continua, chainable and weakly chainable continua, solenoids, circle-like continua, uniformly path connected continua, smooth dendroids, and smooth fans. For smooth fans, we obtain additional characterizations in terms of mappings onto the Cantor and Lelek fans, the containment of these fans, and the number of legs. We also prove that no reciprocating domains exist in the classes of all continua, tree-like continua, planar tree-like continua, and indecomposable continua.
\subsection{Compact metrizable spaces}
In this subsection, we characterize reciprocating domains in the class of
compact metrizable spaces. We prove that a compact metrizable space is a
reciprocating domain precisely when it maps continuously onto the Cantor
set, or equivalently, when it has uncountably many components.

\begin{definition}
	We use $\mathcal{CMS}$ to denote the class of compact metrizable spaces.
\end{definition}

The following well-known theorem and lemma are used to prove Theorem \ref{main:MK}.
\begin{theorem}\label{tanasa}
		Let $X$ be a compact metric space and let $\sim$ be a closed equivalence relation on $X$ such that for each $x\in X$, the equivalence class $[x]$ of $x$ is a closed subset of $X$. Then the quotient space $X/_{\sim}$ is metrizable. 
	\end{theorem}
	\begin{proof}
		See \cite[Theorem 4.2.13]{engelking1}.
	\end{proof} 

\begin{lemma}\label{monique}
Let $X$ be an uncountable totally disconnected compact metric space.
Then there exist clopen sets $U$ and $V$ in $X$ such that $X=U\cup V$, $U\cap V=\emptyset$, and  both $U$ and $V$ are uncountable.
\end{lemma}

\begin{proof}
Since $X$ is a compact metric space, it is second countable. Let $\mathcal{B}$ be a countable basis for its topology, and define
$$
Y=
\left\{
x\in X \mid
\text{every } B\in\mathcal{B} \text{ containing } x
\text{ is uncountable}
\right\}.
$$
For each $x\in X\setminus Y$, $x$ belongs to some countable member of $\mathcal{B}$. Hence
$$
X\setminus Y
\subseteq
\bigcup
\left\{
B\in\mathcal{B} \mid B \text{ is countable}
\right\}.
$$
The union on the right is a countable union of countable sets, so
\(X\setminus Y\) is countable.
It follows that \(Y\) contains at least two distinct points; otherwise
\(X=(X\setminus Y)\cup Y\) would be countable. Choose distinct points
\(x,y\in Y\). Since a totally disconnected compact metric space has a basis of clopen
sets, there exists a clopen set \(B\subseteq X\) such that
$x\in B$ and $y\notin B$. Choose such a set $B$ and let $U=B$ and $V=X\setminus B$. Then both \(U\) and \(V\) are open, and clearly $X=U\cup V$, $U\cap V=\emptyset$, and  both $U$ and $V$ are uncountable.
\end{proof}
\begin{theorem}\label{main:MK}
	Let $X\in \mathcal{CMS}$ and let $C$ be a Cantor set.  The following statements are equivalent:
	\begin{enumerate}
		\item\label{a2} $X$ is a reciprocating domain in $\mathcal{CMS}$.
		\item\label{a1} There is a continuous surjection from $X$ to $C$.
		\item\label{a00} $X$ has uncountably many components.
		\item\label{mu2} $X$ is a universal domain in $\mathcal{CMS}$.
	\end{enumerate}
\end{theorem}

\begin{proof}
		  It follows from Theorem \ref{zelosplosno} that \ref{a2}, \ref{a1} and \ref{mu2} are equivalent. 	The implication from \ref{a00} to \ref{a1} follows from \cite[Theorem 5.1]{KolyadaSemikina2015}.
We next prove that \ref{a1} implies \ref{a00}. Suppose that \(f\colon X\to C\) is a continuous surjection. Since the continuous image of a connected set is connected and every connected subset of \(C\) is a singleton, \(f\) is constant on each component of \(X\). Thus, distinct points of \(C\) must have preimages belonging to distinct components of \(X\). Since \(C\) is uncountable, \(X\) has uncountably many components.
\end{proof}

\subsection{Continua}

In this subsection, we show that there are no reciprocating domains in the
class of all continua. This follows from the nonexistence of universal
continua and the fact that the class of continua is closed under wedges.
\begin{definition}
	A topological space is a \emph{continuum} if it is a non-empty connected compact metrizable space.
	\end{definition}
	\begin{definition}
	We denote by $\mathcal{C}$ the class of all continua.
\end{definition}
We use the following well-known result to prove Theorem \ref{kontinuumi}
\begin{theorem}[Waraszkiewicz]\label{mujika}
		For each $X\in \mathcal{C}$, $X$ is not a universal domain in $\mathcal{C}$.
	\end{theorem}
	\begin{proof}
		The theorem is one of the main results of \cite{W}.
	\end{proof}
\begin{theorem}\label{kontinuumi}
	For each $X\in \mathcal{C}$, $X$ is not a reciprocating domain in $\mathcal{C}$.
\end{theorem}
\begin{proof}
Since  for all $X,Y\in \mathcal{C}$ and for any wedge $W$ of $X$ and $Y$, $W$ is a continuum, it follows from Observation \ref{karin} and Corollary \ref{wedge1} that for any $X\in \mathcal{C}$,
	$$
	X \textup{ is a reciprocating domain in } \mathcal{C} ~~~ \Longleftrightarrow ~~~ X \textup{ is a universal domain in } \mathcal{C}. 
	$$
	 Therefore, it follows from Theorem \ref{mujika} that there are no reciprocating domains in  $\mathcal{C}$.
\end{proof}

\subsection{Peano continua}

In this subsection, we use the Hahn--Mazurkiewicz theorem to characterize
reciprocating domains among Peano continua. We prove that every
nondegenerate Peano continuum is both a reciprocating and a universal
domain in this class.

\begin{definition}
Let $X$ be a continuum. We say that $X$ is a \emph{Peano continuum} if it is a locally connected continuum. \end{definition}
\begin{definition}
We use  $\mathcal{PC}$ to denote the class of all Peano continua.
\end{definition}
The following is a well-known result. 
\begin{theorem}[Hahn-Mazurkiewicz Theorem]\label{hahn}
Let $X$ be a continuum. Then $X$ is a Peano continuum if and only if $X$ is a continuous image of the closed interval $[0,1]$.
\end{theorem}

\begin{proof}
See \cite[Theorems 8.14 and 8.18]{nadler}.
\end{proof}

\begin{corollary}
	$[0,1]$ is a universal domain in $\mathcal{PC}$. 
\end{corollary}
\begin{proof}
	The corollary follows directly from Theorem \ref{hahn}. 
\end{proof}
\begin{observation}\label{bergla}
	Note that for each continuum $X$, if $X$ consists of more than $1$ point, then there is a continuous surjection from $X$ to $[0,1]$. 
\end{observation}
\begin{theorem}
	Let $X\in \mathcal{PC}$. The following statements are equivalent.
	\begin{enumerate}
		\item\label{PC1} $X$ is  a reciprocating domain in $\mathcal{PC}$.
		\item\label{PC2} There is a continuous surjection from $X$ to $[0,1]$.
		\item\label{PC3} $X$ is nondegenerate.
		\item\label{PC4} $X$ is a universal domain in $\mathcal{PC}$.
	\end{enumerate} 
\end{theorem}
\begin{proof}
	By Theorem \ref{zelosplosno}, \ref{PC1}, \ref{PC2} and \ref{PC4} are equivalent. 	Obviously, \ref{PC3} follows from \ref{PC2}. To prove the implication from \ref{PC3} to \ref{PC1}, let $X\in \mathcal{PC}$ be such that $X$ consists of more than $1$ point. By Observation \ref{bergla}, there is a continuous surjection from $X$ to $[0,1]$. Therefore, by Theorem \ref{general:1}, $X$ is a reciprocating domain in $\mathcal{PC}$. 
\end{proof}

\subsection{Chainable and weakly chainable continua}\label{sec:chainable}

In this subsection, we study reciprocating domains in the classes of
chainable and weakly chainable continua. Using the pseudoarc as a
universal domain, we characterize reciprocating domains in terms of the
existence of a continuous surjection onto the pseudoarc. We also examine
the structure of chainable continua that admit such a surjection.
We start the section by defining the class of chainable continua; see Definition~\ref{chainable}. To do that, we need the notion of chains.

\begin{definition}
Let \( X \) be a continuum and let \( U_1, U_2, U_3,\dots, U_n \) be open sets in \( X \). Then we say that the $n$-tuple \( (U_{1}, U_{2}, U_3, \dots, U_{n}) \) is a \emph{chain} in \( X \) if for all \( i,j \in\{1,2,3,\ldots, n\} \), we have
\[
|i-j|\leq 1 ~~~ \Longleftrightarrow ~~~ U_{i} \cap U_{j} \neq \emptyset.
\]
\end{definition}

\begin{definition}
Let \( X \) be a continuum, let $\varepsilon >0$ and let  \( \mathcal U=(U_{1}, U_{2}, U_3, \dots, U_{n}) \) be a chain in \( X \).  We say that $\mathcal U$ is an \emph{$\varepsilon$-chain in $X$}, if for each $i\in\{1,2,3,\ldots,n\}$, the diameter of $U_i$ is less than $\varepsilon$. 
\end{definition}

\begin{definition}
Let \( X \) be a continuum and let  \( \mathcal U=(U_{1}, U_{2}, U_3, \dots, U_{n}) \) be a chain in \( X \).  We say that $\mathcal U$ \emph{covers} $X$, if $\bigcup_{i=1}^{n}U_i=X$. 
\end{definition}

We are now ready to define chainable continua.

\begin{definition}\label{chainable}
Let \( X \) be a continuum. We say that \( X \) is \emph{chainable} if for each $\varepsilon>0$, there is an $\varepsilon$-chain in $X$ that covers $X$. 
\end{definition}

\begin{definition}\label{def:indecomp}
A continuum is \emph{indecomposable} if it is not the union of two proper
subcontinua. It is \emph{hereditarily indecomposable} if every subcontinuum
is indecomposable. A continuum is \emph{nondegenerate} if it contains more
than one point. An \emph{arc} is a continuum homeomorphic to $[0,1]$.
\end{definition}

\begin{definition}\label{mercedes}
A \emph{pseudoarc} is a nondegenerate hereditarily indecomposable chainable
continuum.
\end{definition}

\begin{theorem}[Bing]\label{bing}
Pseudoarcs exist, and any two pseudoarcs are homeomorphic.
Every nondegenerate subcontinuum of a pseudoarc is a pseudoarc.
\end{theorem}
\begin{proof}
Existence and uniqueness follow from \cite{Bing2,Bing3}. A
nondegenerate subcontinuum of a pseudoarc is hereditarily
indecomposable and chainable, so it is a pseudoarc by definition.
\end{proof}

\begin{definition}\label{def:weak-chainable}
A continuum is \emph{weakly chainable} if it is a continuous image of a
chainable continuum.
\end{definition}

\begin{definition}
We use $\mathcal{WC}$  to denote the class of all weakly-chainable continua. 
\end{definition}
\begin{observation}
	Note that if $X\in \mathcal T$ is the degenerate space, then $X\in \mathcal{WC}$.
\end{observation}
The following characterization provides an elegant connection between weakly chainable continua and the pseudoarc.

\begin{theorem}[Lelek]\label{weak}
	Let $X$ be a continuum and let $P$ be a pseudoarc. The following statements are equivalent.
\begin{enumerate}
	\item $X$ is weakly chainable.
	\item There is a continuous surjection from $P$ to $X$.
\end{enumerate}
\end{theorem}

\begin{proof}
	See \cite[Theorem, page 274]{Lelek}.
\end{proof}
\begin{corollary}\label{sedo}
	Let $P$ be a pseudoarc. Then $P$ is a universal domain in $\mathcal{WC}$.
\end{corollary}
\begin{proof}
	The corollary follows directly from Theorem \ref{weak}.
\end{proof}
\begin{theorem}\label{WC}
	Let $X\in \mathcal{WC}$ and let $P$ be a pseudoarc. The following statements are equivalent.
	\begin{enumerate}
		\item\label{WC1} $X$ is a reciprocating domain in $\mathcal{WC}$.
		\item\label{WC2} There is a continuous surjection from $X$ to $P$.
		\item\label{WC3} $X$ is a universal domain in $\mathcal{WC}$. 
	\end{enumerate}
\end{theorem}

\begin{proof}
	By Corollary \ref{sedo}, $P$ is a universal domain in $\mathcal{WC}$. So, the theorem follows from Theorem \ref{zelosplosno}.
\end{proof}

\begin{definition}
	We use  $\mathcal{CC}$ to denote the class of all chainable continua.
\end{definition}

\begin{theorem}\label{jurevkinu}
Let $P$ be a pseudoarc. Then $P$ is a universal as well as a reciprocating domain in $\mathcal{CC}$.
\end{theorem}

\begin{proof}
Since every chainable continuum is also a weakly chainable continuum, it is a continuous image of the pseudoarc $P$. Therefore, $P$ is a universal domain in $\mathcal{CC}$. By Observation~\ref{karin}, it follows that  \( P \) is a reciprocating domain in \( \mathcal{CC} \).
\end{proof}
\begin{theorem}\label{DDJI}
	Let $X\in \mathcal{CC}$ and let $P$ be a pseudoarc. The following statements are equivalent.
	\begin{enumerate}
		\item\label{CC1} $X$ is a reciprocating domain in $\mathcal{CC}$.
		\item\label{CC2} There is a continuous surjection from $X$ to $P$.
		\item\label{CC3} $X$ is a universal domain in $\mathcal{CC}$.
	\end{enumerate}  
\end{theorem}

\begin{proof}
	By Theorem \ref{jurevkinu}, $P$ is a universal as well as a reciprocating domain in $\mathcal{CC}$. So, the theorem follows from Theorem \ref{zelosplosno}. 
\end{proof}

%
Motivated by Theorem \ref{DDJI},  we continue studying in Theorem \ref{bujaha} chainable continua that map onto a pseudoarc. In it, we need a notion of composants of continua as well as the following well-known result. 
\begin{theorem}\label{indecomposable-preimage}
Let \(X\) and \(Y\) be continua, let \(Y\) be indecomposable and let
\(f:X\rightarrow Y\) be a continuous surjection. Then there is an
indecomposable subcontinuum \(K\) of \(X\) such that \(f(K)=Y\).
\end{theorem}

\begin{proof}
The theorem follows directly from \cite[Theorem 4.6]{BlockKeeslingLedis}.
\end{proof}

\begin{definition}
Let \(X\) be a continuum and let \(x\in X\). The \emph{composant} of
\(X\) containing \(x\) is the union of all proper subcontinua of \(X\)
that contain \(x\).
\end{definition}
\begin{observation}
Every nondegenerate indecomposable continuum has infinitely many composants; in fact,
it has uncountably many pairwise disjoint composants. On the other hand,
every decomposable continuum has either exactly one composant or exactly
three composants. For more information on the topic see \cite{nadler}. 
\end{observation}
\begin{theorem}\label{proper-indecomposable-subcontinuum}
Let \(K\) be an indecomposable chainable continuum, let \(P\) be a
pseudoarc, and let
$
f:K\rightarrow P
$
be a continuous surjection. Then \(K\) contains a proper nondegenerate
indecomposable subcontinuum.
\end{theorem}

\begin{proof}
We first prove that there is a proper subcontinuum \(L\) of \(K\) such
that \(f(L)\) is nondegenerate. 
Suppose, towards a contradiction, that \(f(L)\) is degenerate for every
proper subcontinuum \(L\) of \(K\). Let \(C\) be a composant of \(K\),
and fix \(x_0\in C\). For every \(x\in C\), there is a proper
subcontinuum \(L_x\) of \(K\) such that
$x_0,x\in L_x$. For each \(x\in C\), let \(L_x\) be such a continuum. 
By our assumption, \(f(L_x)\) is degenerate. Since it contains both
\(f(x_0)\) and \(f(x)\), it follows that
$f(x)=f(x_0)$. Thus \(f\) is constant on \(C\). Since \(K\) is indecomposable, each composant of \(K\) is dense in
\(K\). Therefore, the continuity of \(f\) implies that \(f\) is
constant on \(K\). This contradicts the assumption that \(f\) is a
surjection onto the nondegenerate continuum \(P\). Hence there is a proper subcontinuum \(L\) of \(K\) such that \(f(L)\)
is nondegenerate. Let $L$ be such a continuum and let
$P_0=f(L)$. Then \(P_0\) is a nondegenerate subcontinuum of the pseudoarc \(P\).
Since every nondegenerate subcontinuum of a pseudoarc is itself a
pseudoarc, \(P_0\) is a pseudoarc. In particular, \(P_0\) is
indecomposable. Consider the restriction
$f|_L:L\rightarrow P_0$.
It is a continuous surjection. By
Theorem~\ref{indecomposable-preimage}, there exists an indecomposable
subcontinuum \(H\) of \(L\) such that
$f(H)=P_0$. Since \(P_0\) is nondegenerate, \(H\) is also nondegenerate. Moreover,
$H\subseteq L\subseteq K$ and $L\neq K$. Therefore \(H\) is a proper nondegenerate indecomposable subcontinuum
of \(K\).
\end{proof}
\begin{corollary}\label{decreasing-indecomposable-sequence}
Let \(X\) be a chainable continuum, let \(P\) be a pseudoarc, and let
$f:X\rightarrow P$ be a continuous surjection. Then there exists a sequence $(K_n)$ 
of nondegenerate indecomposable subcontinua of \(X\) such that for each positive integer $n$, $
K_{n+1}\subseteq K_n$ and $K_{n+1}\neq K_n$. 
\end{corollary}

\begin{proof}
By Theorem~\ref{indecomposable-preimage}, there exists an
indecomposable subcontinuum \(K_1\) of \(X\) such that
$f(K_1)=P$. Since \(P\) is nondegenerate, \(K_1\) is nondegenerate. Moreover, since
\(X\) is chainable, \(K_1\) is also chainable. We construct the sequence \((K_n)\) inductively. Suppose that $n$ is a positive integer and that \(K_n\) has already been constructed in such a way that \(K_n\) is
a nondegenerate indecomposable chainable continuum and $f(K_n)$ is a pseudoarc. Let $P_n=f(K_n)$.
Consider the continuous surjection
\[
f|_{K_n}:K_n\rightarrow P_n.
\]
By Theorem~\ref{proper-indecomposable-subcontinuum} and its proof, there
exists a proper nondegenerate indecomposable subcontinuum
\(K_{n+1}\) of \(K_n\) such that
$f(K_{n+1})$ 
is a nondegenerate subcontinuum of \(P_n\). Let $K_{n+1}$ be such a continuum. Since \(K_{n+1}\) is a subcontinuum of the chainable continuum \(K_n\),
it is chainable. Moreover, every nondegenerate subcontinuum of a
pseudoarc is itself a pseudoarc. Therefore,
$f(K_{n+1})$
is a pseudoarc. It follows by induction that there exists a sequence
$(K_n)$ of nondegenerate indecomposable subcontinua of \(X\) such that for each positive integer $n$, $
K_{n+1}\subseteq K_n$ and $K_{n+1}\neq K_n$. 
\end{proof}

\begin{theorem}\label{bujaha}
Let $X$ be a chainable continuum, let $P$ be a pseudoarc, and let $f\colon X\rightarrow P$ be a continuous surjection. For every composant $C$ of $P$, define
$$
\mathcal{K}_C
=
\left\{
Y\subseteq X \mid
Y \text{ is a subcontinuum of }X
\text{ and } f(Y)\subseteq C
\right\}.
$$
Then, for every composant $C$ of $P$, there exists a sequence $(K_n)$ 
of nondegenerate indecomposable continua in $\mathcal{K}_C$ such that for each positive integer $n$, $
K_{n+1}\subseteq K_n$ and $K_{n+1}\neq K_n$.
\end{theorem}
\begin{proof}
Let $C$ be a composant of $P$. Choose $p\in C$ and
$x\in f^{-1}(\{p\})$. By the order-arc theorem for the hyperspace $C(X)$ of non-empty
subcontinua of $X$ (see \cite{nadler}), there exists a continuous map
$\alpha\colon[0,1]\to C(X)$, where $C(X)$ is endowed
with the Hausdorff metric, such that
\[
\alpha(0)=\{x\} ~  \text{ and } ~ \alpha(1)=X,
\]
and $\alpha(s)\subseteq\alpha(t)$ whenever $s\leq t$.
The function from $[0,1]$ to $[0,\diam(P)]$, defined by 
\[
t\longmapsto \operatorname{diam} f(\alpha(t))
\]
for each $t\in [0,1]$, is continuous, takes the value $0$ at $t=0$, and takes
the value $\operatorname{diam}P>0$ at $t=1$.
Consequently, there exists $t_0\in(0,1)$ such that
\[
0<\operatorname{diam} f(\alpha(t_0))
  <\operatorname{diam}P.
\]
Set $H=\alpha(t_0)$ and $P_0=f(H)$.
Then $P_0$ is a proper nondegenerate subcontinuum
of $P$ containing $p$. Thus $P_0\subseteq C$, and
$P_0$ is a pseudoarc.
By Theorem \ref{indecomposable-preimage}, there exists
an indecomposable subcontinuum $K\subseteq H$ such that
$f(K)=P_0$. In particular, $K$ is nondegenerate and
chainable. Apply Corollary
\ref{decreasing-indecomposable-sequence} to
$f|_K\colon K\to P_0$ to obtain a sequence $(K_n)$
of nondegenerate indecomposable subcontinua of $K$
such that
\[
K_{n+1}\subseteq K_n ~  \text{ and } ~ K_{n+1}\neq K_n
\]
for every positive integer $n$.
Since $f(K_n)\subseteq P_0\subseteq C$, every $K_n$
belongs to $\mathcal K_C$.
\end{proof}

\subsubsection{Chainable domains without pseudoarcs}\label{sec:pa-construction-start}

Theorem~\ref{bujaha} supplies many indecomposable subcontinua in a
chainable domain mapping onto a pseudoarc. We now show that none of
these subcontinua need be hereditarily indecomposable. By
Theorem~\ref{DDJI}, the construction also gives a reciprocating and
universal domain in $\mathcal{CC}$ containing no pseudoarc.
The definitions and lemmas below give the complete construction.

\begin{definition}
Let $(X_n)$ be a sequence of continua and, for each positive integer
$n$, let $f_n:X_{n+1}\to X_n$ be a continuous map. We call
$(X_n,f_n)$ an \emph{inverse sequence}. We define its
\emph{inverse limit} by
\[
\invlim(X_n,f_n)=
\left\{(x_n)\in\prod_{n=1}^{\infty}X_n\ \middle|\
f_n(x_{n+1})=x_n\text{ for each positive integer }n\right\},
\]
with the topology inherited from $\prod_{n=1}^{\infty}X_n$.
For each positive integer $n$, the map
$\pi_n:\invlim(X_n,f_n)\to X_n$, defined by
$\pi_n((x_i))=x_n$ for each $(x_i)\in\invlim(X_n,f_n)$,
is called the $n$th coordinate projection.
\end{definition}

\begin{definition}
Let $X,Y$ be continua. A continuous surjection $f:X\to Y$ is
\emph{monotone} if each fiber $f^{-1}(\{y\})$, $y\in Y$, is connected.
\end{definition}
\begin{definition}
	Let $\id:[0,1]\to[0,1]$ be the identity map. We define
$R,T:[0,1]\to[0,1]$ by $R(t)=1-t$ and 
\[
T(t)=\begin{cases}2t,&0\le t\le 1/2,\\2-2t,&1/2\le t\le1.\end{cases}
\]
\end{definition}

\begin{definition}
Let $a,b\in\mathbb R$ be such that $a<b$, and let $E=[a,b]$.
The map $\beta:[0,1]\to E$ defined by
\[
\beta(t)=a+(b-a)t
\]
for each $t\in[0,1]$ is called the \emph{increasing affine
homeomorphism from $[0,1]$ onto $E$}.
\end{definition}

\begin{definition}
A topological space $X$ is said to be \emph{zero-dimensional}
if it has a basis consisting of clopen sets. Equivalently, for
each $x\in X$ and each open set $U\subseteq X$ containing $x$,
there is a clopen set $V\subseteq X$ such that
$x\in V\subseteq U$.
\end{definition}

\begin{definition}
Let $a,b\in\mathbb R$ be such that $a<b$, and let
$g:[a,b]\to\mathbb R$ be a continuous map. We say that $g$ is
\emph{piecewise linear} if there are a positive integer $m$ and
points $a=c_0<c_1<\cdots<c_m=b$ such that $g$ is linear on
$[c_{i-1},c_i]$ for each $i\in\{1,2,\ldots,m\}$.
We say that $g$ is \emph{nowhere constant} if it is not constant
on any nondegenerate interval.
\end{definition}

\begin{definition}
Let $a,b\in\mathbb R$ be such that $a<b$, and let
$g:[a,b]\to\mathbb R$ be a continuous piecewise linear map
that is not constant on any nondegenerate interval.
A point $t\in(a,b)$ is called a \emph{turning point} of $g$
if $g$ has a local maximum or a local minimum at $t$.
\end{definition}

\begin{definition}
Let $a,b\in\mathbb R$ be such that $a<b$, and let
$f:[a,b]\to\mathbb R$. We say that $f$ is
\emph{right-continuous at $t\in[a,b)$} if
\[
\lim_{s\to t^+}f(s)=f(t).
\]
Equivalently, for each $\varepsilon>0$, there is a $\delta>0$
such that $|f(s)-f(t)|<\varepsilon$ for each $s\in[a,b]$
satisfying $t\le s<t+\delta$.
We say that $f$ is \emph{right-continuous} if it is
right-continuous at each $t\in[a,b)$.
\end{definition}

\subsubsection{Expanding the coordinate intervals}\label{pa-sec:construction}

The following lemma, which is used in the proof of
Lemma~\ref{pa-lem:arc-extraction}, is a special case of a more general result of Read
\cite{Read1974}, which states that every arc-like continuum belongs
to Class~W, that is, every continuous surjection from a continuum
onto an arc-like continuum is weakly confluent. This result is also
recorded in \cite[Theorem~1.2]{GrispolakisTymchatyn1978}.
Since every interval is arc-like, weak confluence implies that
every closed subinterval of the range is the image of some
subcontinuum of the domain. We include a direct proof for completeness.

\begin{lemma}\label{pa-lem:crossing}
Let $H$ be a continuum, let $c,d\in\mathbb R$ be such that $c<d$,
and let $q:H\to[c,d]$ be a continuous surjection. Also, let
$a,b\in[c,d]$ be such that $a<b$. Then there is a subcontinuum
$L$ of $H$ such that $q(L)=[a,b]$.
\end{lemma}
\begin{proof}
Let $Z=q^{-1}([a,b])$, let $E=q^{-1}(\{a\})$, and let
$F=q^{-1}(\{b\})$. We prove that some component of $Z$ meets both
$E$ and $F$. Suppose, towards a contradiction, that no component
of $Z$ meets both $E$ and $F$.

Let $e\in E$. Since $Z$ is a compact Hausdorff space, the component of $Z$
containing $e$ is the intersection of all clopen subsets of $Z$
that contain $e$; see \cite[Lemma 5.12.11]{StacksComponents}. By our assumption, this component is
disjoint from $F$. Therefore, for each $z\in F$, there is a
clopen subset $U_{e,z}$ of $Z$ such that $e\in U_{e,z}$ and
$z\notin U_{e,z}$. Since $F$ is compact, there are finitely many
such sets whose intersection is disjoint from $F$. Let $U_e$
be their intersection. Then $U_e$ is clopen in $Z$,
$e\in U_e$, and $U_e\cap F=\emptyset$.

For each $e\in E$, choose such a clopen set $U_e$. Then
$\{U_e\mid e\in E\}$ is an open cover of $E$. Since $E$ is compact,
finitely many of the sets $U_e$, $e\in E$, cover $E$. Let $U$ be their union. Then $U$ is clopen in $Z$,
$E\subseteq U$, and $U\cap F=\emptyset$. Note that
\[
q^{-1}([c,a])\cup U
\quad\text{and}\quad
(Z\setminus U)\cup q^{-1}([b,d])
\]
are disjoint non-empty closed subsets of $H$ whose union is $H$.
This contradicts the assumption that $H$ is connected. Therefore,
there is a component $L$ of $Z$ that meets both $E$ and $F$.
Let $L$ be such a component. Since $Z$ is compact, $L$ is a
subcontinuum of $H$. Also, $q(L)\subseteq[a,b]$ and
$a,b\in q(L)$. Since $q(L)$ is connected, it follows that
$q(L)=[a,b]$.
\end{proof}

The following lemma is used to study the subcontinua of the fibers
of the map constructed in Lemma~\ref{pa-lem:limit-fibers}.

\begin{lemma}\label{pa-lem:arcs}
Let $(E_n,h_n)$ be an inverse sequence such that for each positive
integer $n$, $E_n$ is a nondegenerate compact interval and
$h_n:E_{n+1}\to E_n$ is a continuous map. For each positive integer
$n$, let $b_n:[0,1]\to E_n$ be the increasing affine homeomorphism,
and let
\[
\varphi_n=b_n^{-1}\circ h_n\circ b_{n+1}.
\]
If for each positive integer $n$,
$\varphi_n\in\{\id,R,T,R\circ T\}$, then every nondegenerate
subcontinuum of $\invlim(E_n,h_n)$ contains an arc.
\end{lemma}

\begin{proof}
Let $Z=\invlim(E_n,h_n)$ and let
$W=\invlim([0,1],\varphi_n)$. First, we prove that $Z$ and $W$
are homeomorphic. Let $\Theta:W\to Z$ be defined by
\[
\Theta((t_n))=(b_n(t_n))
\]
for each $(t_n)\in W$. Note that
$h_n\circ b_{n+1}=b_n\circ\varphi_n$ for each positive integer
$n$. Therefore, $\Theta$ is well-defined and continuous. Its inverse
is the continuous map defined by
\[
(t_n)\longmapsto(b_n^{-1}(t_n))
\]
for each $(t_n)\in Z$. Hence, $\Theta$ is a homeomorphism.
It follows that it suffices to prove that every nondegenerate
subcontinuum of $W$ contains an arc.

Let $H$ be a nondegenerate subcontinuum of $W$ and, for each
positive integer $n$, let $\pi_n:W\to[0,1]$ be the $n$th
coordinate projection and let $A_n=\pi_n(H)$. Then each $A_n$
is a compact interval and $\varphi_n(A_{n+1})=A_n$.
Since $H$ is nondegenerate, there is a positive integer $m$
such that $A_m$ is nondegenerate. Let $m$ be such an integer
and let $B_m\subseteq A_m\cap(0,1)$ be a nondegenerate compact
interval.

We construct intervals $B_n$, $n\ge m$, inductively. Suppose that
$n\ge m$ and that a nondegenerate compact interval
$B_n\subseteq A_n\cap(0,1)$ has already been constructed.
We prove that there is a nondegenerate compact interval
$B_{n+1}\subseteq A_{n+1}\cap(0,1)$ such that
$\varphi_n|_{B_{n+1}}$ is a homeomorphism from $B_{n+1}$ onto $B_n$.
If $\varphi_n\in\{\id,R\}$, let $B_{n+1}=\varphi_n^{-1}(B_n)$.
Then $B_{n+1}$ has the required properties.

Suppose that $\varphi_n=T$. If $A_{n+1}$ is contained in
$[0,1/2]$ or in $[1/2,1]$, then $T|_{A_{n+1}}$ is a
homeomorphism onto $A_n$. In this case, let $B_{n+1}$ be the
inverse image of $B_n$ under this homeomorphism. Otherwise,
$1/2\in A_{n+1}$ and at least one of the intervals
\[
A_{n+1}\cap[0,1/2]
\quad\text{and}\quad
A_{n+1}\cap[1/2,1]
\]
is mapped by $T$ onto $A_n$. Let $C$ be such an interval, and let
$B_{n+1}=(T|_C)^{-1}(B_n)$. Again, $B_{n+1}$ has the required
properties. In both cases, $B_{n+1}\subseteq(0,1)$ because
$B_n\subseteq(0,1)$. The case $\varphi_n=R\circ T$ follows in
the same way. Figure~\ref{pa-fig:fourmaps} illustrates this construction.
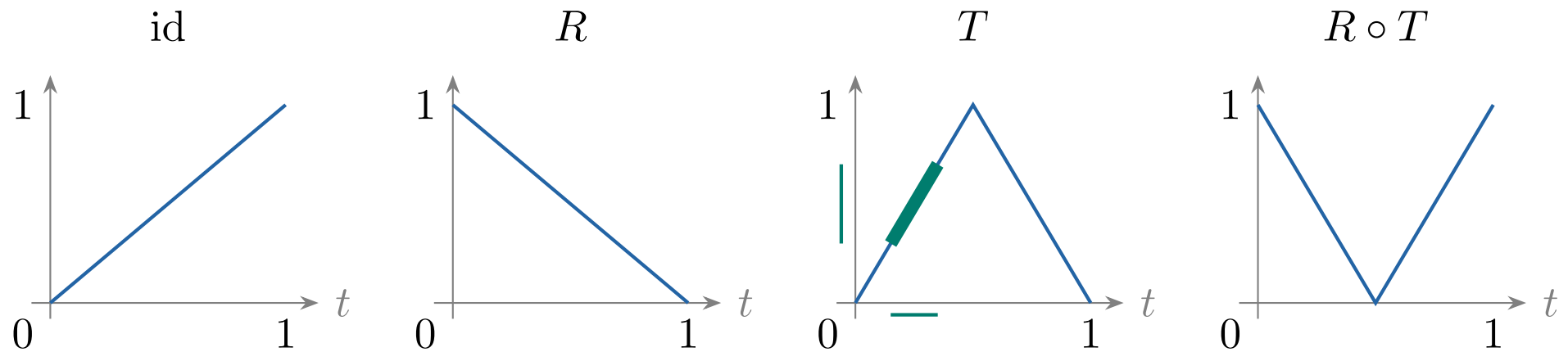
\begin{figure}[!ht]
\centering
\begin{tikzpicture}[x=1.9cm,y=1.6cm,font=\small]
\foreach \k/\lab in {0/{$\id$},1/{$R$},2/{$T$},3/{$R\circ T$}} {
  \begin{scope}[xshift=\k*3.25cm]
  \draw[axis] (-.08,0)--(1.14,0) node[right] {$t$};
  \draw[axis] (0,-.08)--(0,1.15);
  \node[below left] at (0,0) {$0$};
  \node[below] at (1,0) {$1$};
  \node[left] at (0,1) {$1$};
  \node at (.5,1.4) {\lab};
  \end{scope}
}
\draw[graphline] (0,0)--(1,1);
\begin{scope}[xshift=3.25cm]\draw[graphline] (0,1)--(1,0);\end{scope}
\begin{scope}[xshift=6.5cm]
  \draw[graphline] (0,0)--(.5,1)--(1,0);
  \draw[branchgreen,line width=3pt] (.15,.3)--(.35,.7);
  \draw[branchgreen,thick] (-.06,.3)--(-.06,.7);
  \draw[branchgreen,thick] (.15,-.06)--(.35,-.06);
\end{scope}
\begin{scope}[xshift=9.75cm]\draw[graphline] (0,1)--(.5,0)--(1,1);\end{scope}
\end{tikzpicture}
\caption{The maps $\id$, $R$, $T$, and $R\circ T$. The marked branch
of $T$ maps $B_{n+1}$ homeomorphically onto $B_n$.}
\label{pa-fig:fourmaps}
\end{figure}
It follows by induction that the intervals $B_n$, $n\ge m$, have
been constructed. Since all the restricted bonding maps are
homeomorphisms, their inverse limit is an arc. Define its coordinates
with index less than $m$ using the original bonding maps. We obtain
an arc whose $n$th coordinate belongs to $A_n$ for every positive
integer $n$. We prove that this arc is contained in $H$.

Let $(t_n)$ be a point of this arc. For each positive integer $n$, let
\[
H_n=\{z\in H:\pi_n(z)=t_n\}.
\]
Since $t_n\in A_n$, each $H_n$ is non-empty and compact. Moreover,
$H_{n+1}\subseteq H_n$ for each positive integer $n$, by the
compatibility equations for the coordinates. Therefore,
$\bigcap_{n\ge1}H_n\ne\emptyset$. Every point in this intersection
has coordinates $(t_n)$, so $(t_n)\in H$. Hence, the arc is
contained in $H$. Therefore, every nondegenerate subcontinuum of
$W$, and consequently of $Z$, contains an arc.
\end{proof}

In the rest of this section, we use the following notation. Let $(I_n,g_n)$
be an inverse sequence such that for each positive integer $n$,
$I_n=[0,1]$ and $g_n:I_{n+1}\to I_n$ is a continuous surjection
that is piecewise linear and nowhere constant. Let
\[
Y=\invlim(I_n,g_n)
\]
and, for each positive integer $n$, let $\pi_n:Y\to I_n$ be the
$n$th coordinate projection.

We first consider the fibers of the coordinate projections of $Y$.

\begin{lemma}\label{pa-lem:coordinate-fibers}
Let $n$ be a positive integer and let $t\in I_n$. Then
$\pi_n^{-1}(\{t\})$ is a zero-dimensional compact metric space.
\end{lemma}
\begin{proof}
Let $n$ be a positive integer and let $t\in I_n$. For each integer
$m>n$, let
\[
S_m=(g_n\circ g_{n+1}\circ\cdots\circ g_{m-1})^{-1}(\{t\}).
\]
Since each $g_i$ is finite-to-one, every $S_m$ is finite.
For each $x\in\pi_n^{-1}(\{t\})$ and each $m>n$, the $m$th
coordinate of $x$ belongs to $S_m$. Also, the coordinates of $x$
with index less than $n$ are determined by $t$. Therefore,
$\pi_n^{-1}(\{t\})$ is homeomorphic to a subspace of
$\prod_{m>n}S_m$. Since each $S_m$ is finite and discrete,
this product is zero-dimensional. It follows that
$\pi_n^{-1}(\{t\})$ is zero-dimensional. It is compact and
metrizable because it is a closed subset of $Y$.
\end{proof}

To define the new coordinate intervals, we use the following lemma.

\begin{lemma}\label{pa-lem:dense-sets}
There is a sequence $(D_n)$ such that for each positive integer $n$,
$D_n$ is a countable dense subset of $I_n$, $0,1\in D_n$,
$D_{n+1}$ contains every turning point of $g_n$, and
$g_n^{-1}(D_n)=D_{n+1}$.
\end{lemma}
\begin{proof}
For each positive integer $n$, let $E_n$ be a countable dense subset
of $I_n$ containing $0,1$ and, if $n>1$, all turning points of
$g_{n-1}$. We construct sets $D_n^{(k)}$ for all positive integers
$n$ and all nonnegative integers $k$ inductively. For each positive
integer $n$, let $D_n^{(0)}=E_n$. Suppose that the sets $D_n^{(k)}$
have already been constructed for every positive integer $n$.
We define
\[
D_n^{(k+1)}=D_n^{(k)}\cup g_n(D_{n+1}^{(k)})
\cup\begin{cases}
\emptyset,&n=1,\\
g_{n-1}^{-1}(D_{n-1}^{(k)}),&n>1.
\end{cases}
\]
Since each $g_n$ is finite-to-one, all the sets $D_n^{(k)}$ are
countable. For each positive integer $n$, let
\[
D_n=\bigcup_{k\ge0}D_n^{(k)}.
\]
Then $D_n$ is countable. Since $E_n\subseteq D_n$, it follows that
$D_n$ is dense in $I_n$, contains $0,1$, and contains all turning
points of $g_{n-1}$ if $n>1$.

Let $n$ be a positive integer. By the construction,
$g_n(D_{n+1})\subseteq D_n$ and $g_n^{-1}(D_n)\subseteq D_{n+1}$.
The first inclusion also implies that
$D_{n+1}\subseteq g_n^{-1}(D_n)$. Therefore,
\begin{equation}\label{pa-eq:invariance}
g_n^{-1}(D_n)=D_{n+1}.
\end{equation}
Hence, the sequence $(D_n)$ has the required properties.
\end{proof}

Next, we replace each point of $D_n$ by a nondegenerate interval.

\begin{lemma}\label{pa-lem:expansion}
Let $(D_n)$ be a sequence as in Lemma~\ref{pa-lem:dense-sets} and,
for each positive integer $n$, let $J_n=[0,2]$. Then for each
positive integer $n$, there is a continuous monotone surjection
$r_n:J_n\to I_n$ such that for every $t\in I_n$,
$F_{n,t}=r_n^{-1}(\{t\})$ is a (possibly degenerate) compact interval and
\[
F_{n,t}\text{ is nondegenerate}\quad\Longleftrightarrow\quad t\in D_n.
\]
\end{lemma}
\begin{proof}
Let $n$ be a positive integer and let
$D_n=\{d_{n,1},d_{n,2},d_{n,3},\ldots\}$ be an enumeration
without repetitions. For each $t\in I_n$, we define
\[
\ell_n(t)=t+\sum_{d_{n,k}<t}2^{-k},
~~~
u_n(t)=t+\sum_{d_{n,k}\le t}2^{-k}
 ~~~ \text{ and } ~~~
F_{n,t}=[\ell_n(t),u_n(t)].
\]
We prove that the intervals $F_{n,t}$, $t\in I_n$, form a
partition of $J_n$. Let $s,t\in I_n$ be such that $s<t$.
Then
\[
\ell_n(t)-u_n(s)
=t-s+\sum_{s<d_{n,k}<t}2^{-k}>0.
\]
Therefore, $F_{n,s}\cap F_{n,t}=\emptyset$.
Moreover, $\ell_n(0)=0$, $u_n(1)=2$, and the increasing function
$t\mapsto u_n(t)$ is right-continuous, with left-hand limit
$\ell_n(t)$ at every $t>0$. These facts follow directly from
convergence of the series $\sum_{k\ge1}2^{-k}$. It follows that
the intervals $[\ell_n(t),u_n(t)]$, $t\in I_n$, cover $[0,2]$.
Hence, they form a partition of $J_n$.

We define $r_n:J_n\to I_n$ by $r_n(x)=t$ whenever $x\in F_{n,t}$.
By the preceding paragraph, $r_n$ is well-defined and surjective.
To prove continuity, let $x,y\in J_n$ be such that $x<y$, and
let $s=r_n(x)$ and $t=r_n(y)$. Then $s\le t$. If $s<t$, we have
\[
t-s\le\ell_n(t)-u_n(s)\le y-x.
\]
The inequality $|r_n(x)-r_n(y)|\le|x-y|$ also holds if $s=t$.
Therefore, $r_n$ is continuous. Its fibers are the intervals
$F_{n,t}$, so $r_n$ is monotone.

If $t=d_{n,k}\in D_n$, then $u_n(t)-\ell_n(t)=2^{-k}>0$.
If $t\notin D_n$, then $u_n(t)=\ell_n(t)$. Thus $F_{n,t}$ is
nondegenerate if and only if $t\in D_n$, as required.
Figure~\ref{pa-fig:expansion} shows the map $r_n$.

\begin{figure}[!ht]
\centering
\begin{tikzpicture}[x=1cm,y=1cm,font=\small]
\draw[thick] (0,2)--(12,2);
\draw[fiber] (1.5,2)--(3,2);
\draw[fiber] (5.25,2)--(6,2);
\draw[fiber] (8.25,2)--(9.35,2);
\node[left] at (0,2) {$J_n$};
\node[above] at (2.25,2.1) {$F_{n,t_1}$};
\node[above] at (5.625,2.1) {$F_{n,t_2}$};
\node[above] at (8.8,2.1) {$F_{n,t_3}$};
\draw[thick] (0,0)--(12,0);
\node[left] at (0,0) {$I_n$};
\foreach \a/\b/\c/\lab in {1.5/3/2.25/1,5.25/6/5.625/2,8.25/9.35/8.8/3} {
  \draw[gray,dashed] (\a,1.9)--(\c,.12);
  \draw[gray,dashed] (\b,1.9)--(\c,.12);
  \fill[fiberblue] (\c,0) circle (2pt);
  \node[below] at (\c,0) {$t_{\lab}\in D_n$};
}
\draw[maparrow] (11.2,1.65)--node[right] {$r_n$} (11.2,.35);
\fill (10.1,2) circle (1.5pt);
\fill (10.1,0) circle (1.5pt);
\draw[gray,dashed] (10.1,1.9)--(10.1,.12);
\node[below] at (10.1,-.52) {$t\notin D_n$};
\end{tikzpicture}
\caption{The map $r_n$. Only three nondegenerate fibers are shown.}
\label{pa-fig:expansion}
\end{figure}
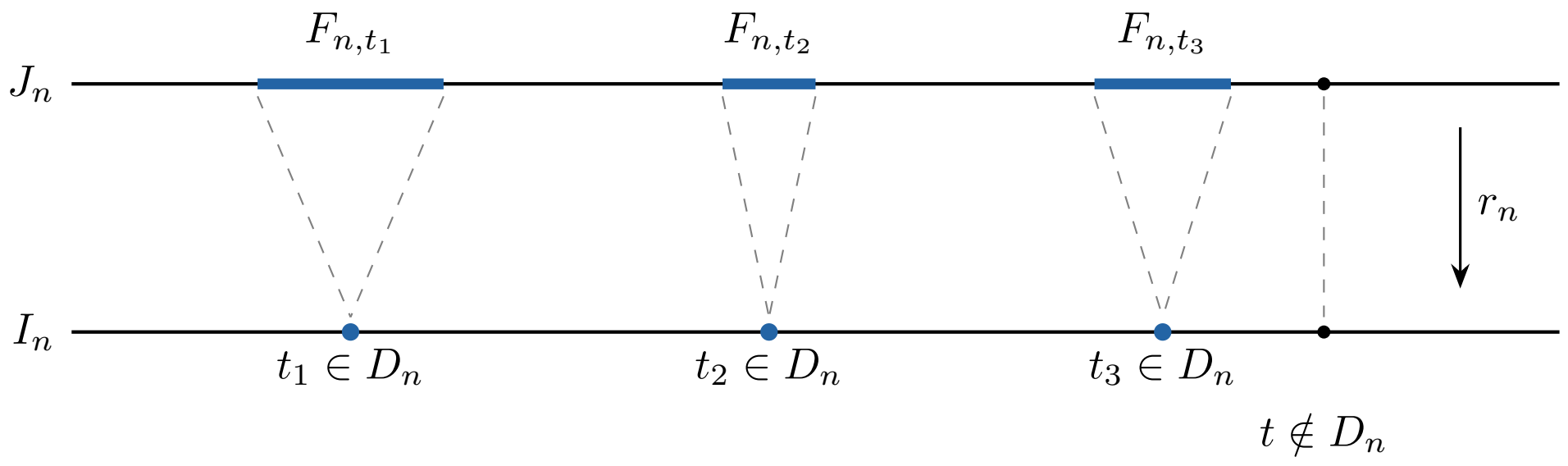
\end{proof}
\begin{definition}\label{pa-def:fiber-affine}
Let $(D_n)$ be as in Lemma~\ref{pa-lem:dense-sets}, and let
$F_{n,t}=[\ell_n(t),u_n(t)]$ be the intervals constructed in
Lemma~\ref{pa-lem:expansion}. For each positive integer $n$ and
each $t\in D_n$, we define
\[
a_{n,t}:[0,1]\to F_{n,t}
\]
to be the increasing affine homeomorphism defined by
\[
a_{n,t}(s)=\ell_n(t)
          +\bigl(u_n(t)-\ell_n(t)\bigr)s
\]
for each $s\in[0,1]$.
\end{definition}
We now define the bonding maps on the new coordinate intervals.

\begin{lemma}\label{pa-lem:lifting}
Let $(D_n)$ be as in Lemma~\ref{pa-lem:dense-sets}, and for each
positive integer $n$, let $r_n$, $J_n$, and the family
$(F_{n,s})_{s\in I_n}$ be as in Lemma~\ref{pa-lem:expansion}.
Also, let $a_{n,s}$, $s\in D_n$, be as in
Definition~\ref{pa-def:fiber-affine}.
Then there is a sequence $(G_n)$ of continuous surjections
$G_n:J_{n+1}\to J_n$ such that for each positive integer $n$,
\[
r_n\circ G_n=g_n\circ r_{n+1}
\quad\text{and}\quad
G_n(F_{n+1,s})=F_{n,g_n(s)}
\text{ for each }s\in I_{n+1}.
\]
Also, for each positive integer $n$ and each $s\in D_{n+1}$,
the map
\[
\psi_{n,s}=a_{n,g_n(s)}^{-1}\circ
\bigl(G_n|_{F_{n+1,s}}\bigr)\circ a_{n+1,s}
\]
belongs to $\{\id,R,T,R\circ T\}$.
\end{lemma}
\begin{proof}
Let $n$ be a positive integer. We define a map
$G_n:J_{n+1}\to J_n$ such that
\begin{equation}\label{pa-eq:commute}
r_n\circ G_n=g_n\circ r_{n+1}.
\end{equation}
If $s\notin D_{n+1}$, then $g_n(s)\notin D_n$ by
\eqref{pa-eq:invariance}. In this case, define $G_n$ to send
the unique point of $F_{n+1,s}$ to the unique point of
$F_{n,g_n(s)}$.
If $s\in D_{n+1}$, define
\[
G_n|_{F_{n+1,s}}
=a_{n,g_n(s)}\circ\psi_{n,s}\circ a_{n+1,s}^{-1},
\]
where $\psi_{n,s}$ is chosen as follows for $s\in(0,1)$:
\[
\psi_{n,s}=
\begin{cases}
\id,&g_n\text{ is strictly increasing near }s,\\
R,&g_n\text{ is strictly decreasing near }s,\\
T,&s\text{ is a local maximum of }g_n,\\
R\circ T,&s\text{ is a local minimum of }g_n.
\end{cases}
\]
Here, increasing or decreasing near $s$ means that the restriction
of $g_n$ to some neighborhood of $s$ is strictly increasing or strictly
decreasing, respectively. At $s=0$ or $s=1$, use $\id$ or $R$ according to the one-sided
monotonicity of $g_n$. These cases exhaust all possibilities,
since $g_n$ has finitely many nonconstant linear pieces. Figure~\ref{pa-fig:lifting} illustrates the four cases.

\begin{figure}[p]
\centering
\begingroup
\newcommand{\liftingpanel}[5]{%
  \node[font=\normalsize] at (6,3.25) {#1};

  \begin{scope}[x=3.5cm,y=2.3cm]
    \draw[axis] (0,0)--(1.08,0);
    \draw[axis] (0,0)--(0,1.08);
    \draw[graphline] #2;
    \draw[gray,dashed] (.5,0)--(.5,#5)--(0,#5);
    \fill[fiberblue] (.5,#5) circle (1.5pt);
    \node[below] at (.5,0) {$s$};
    \node[left] at (0,#5) {$g_n(s)$};
    \node[above] at (.5,1.08) {$g_n$};
  \end{scope}

  \draw[maparrow] (4.8,1.15)--(6.1,1.15);

  \begin{scope}[xshift=8cm,x=3.5cm,y=2.3cm]
    \fill[fiberblue!9] (.3,.2) rectangle (.7,.8);
    \draw[axis] (0,0)--(1.08,0);
    \draw[axis] (0,0)--(0,1.08);
    \draw[graphline] #3;
    \draw[gray,dashed] (.3,0)--(.3,.8);
    \draw[gray,dashed] (.7,0)--(.7,.8);
    \draw[gray,dashed] (0,.2)--(.7,.2);
    \draw[gray,dashed] (0,.8)--(.7,.8);
    \draw[fiber] (.3,0)--(.7,0);
    \draw[fiber] (0,.2)--(0,.8);
    \node[below] at (.5,0) {$F_{n+1,s}$};
    \node[left,font=\small] at (0,.2)
      {$\ell_n(g_n(s))$};
    \node[left,font=\small] at (0,.8)
      {$u_n(g_n(s))$};
    \node[above] at (.5,1.08) {$G_n$};
  \end{scope}

  \node[font=\normalsize] at (6,-.65)
    {$\psi_{n,s}=#4$};
}

\begin{tikzpicture}[font=\small]
  \begin{scope}
    \liftingpanel
      {(a) Increasing near $s$}
      {(.08,.08)--(.92,.92)}
      {(.05,.05)--(.3,.2)--(.7,.8)--(.95,.95)}
      {\id}
      {.5}
  \end{scope}

  \begin{scope}[yshift=-4.4cm]
    \liftingpanel
      {(b) Decreasing near $s$}
      {(.08,.92)--(.92,.08)}
      {(.05,.95)--(.3,.8)--(.7,.2)--(.95,.05)}
      {R}
      {.5}
  \end{scope}

  \begin{scope}[yshift=-8.8cm]
    \liftingpanel
      {(c) Local maximum at $s$}
      {(.08,.15)--(.5,.8)--(.92,.15)}
      {(.05,.05)--(.3,.2)--(.5,.8)--(.7,.2)--(.95,.05)}
      {T}
      {.8}
  \end{scope}

  \begin{scope}[yshift=-13.2cm]
    \liftingpanel
      {(d) Local minimum at $s$}
      {(.08,.85)--(.5,.2)--(.92,.85)}
      {(.05,.95)--(.3,.8)--(.5,.2)--(.7,.8)--(.95,.95)}
      {R\circ T}
      {.2}
  \end{scope}
\end{tikzpicture}
\endgroup
\caption{The four cases in the construction of $G_n$: (a) $\id$;
(b) $R$; (c) $T$; (d) $R\circ T$.
The shaded rectangles correspond to
$F_{n+1,s}\times F_{n,g_n(s)}$.
At $s=0$ or $s=1$, the corresponding one-sided increasing
or decreasing case is used.}
\label{pa-fig:lifting}
\end{figure}
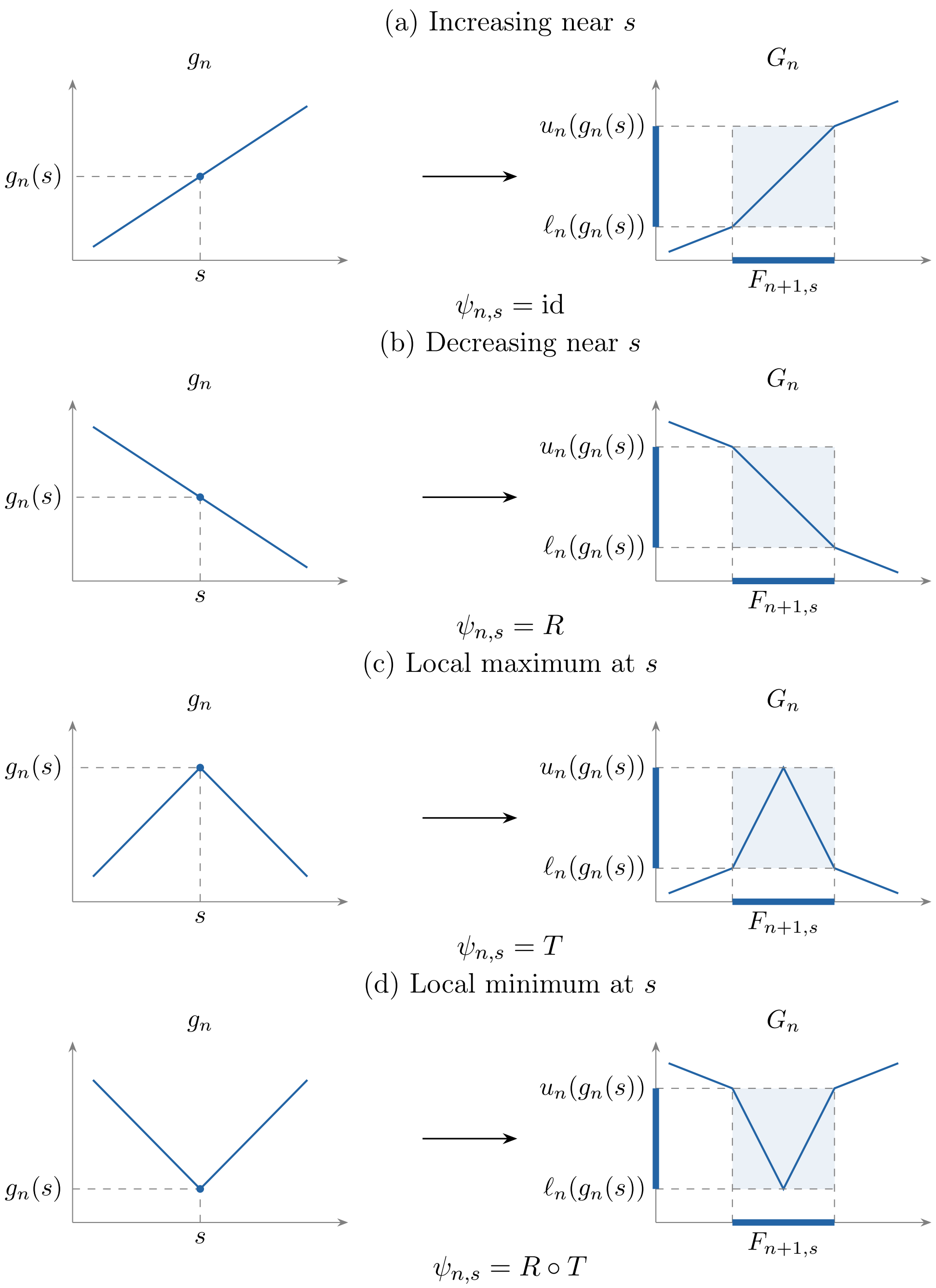

We prove that $G_n$ is continuous. Let $t\in J_{n+1}$, let
$(t_j)$ be a sequence in $J_{n+1}$ such that
$\lim_{j\to\infty}t_j=t$, and let $s=r_{n+1}(t)$.
By \eqref{pa-eq:commute}, continuity of $r_n,r_{n+1},g_n$,
and compactness of $J_n$, every cluster point of
$(G_n(t_j))$ belongs to $F_{n,g_n(s)}$.

First, suppose that $F_{n+1,s}$ is a singleton. Then
$s\notin D_{n+1}$. Since $g_n^{-1}(D_n)=D_{n+1}$ by
Lemma~\ref{pa-lem:dense-sets}, it follows that $g_n(s)\notin D_n$.
Therefore, $F_{n,g_n(s)}$ is also a singleton.
Since every cluster point of $(G_n(t_j))$ belongs to
$F_{n,g_n(s)}=\{G_n(t)\}$ and $J_n$ is compact, it follows
that $\lim_{j\to\infty}G_n(t_j)=G_n(t)$.

Next, suppose that $F_{n+1,s}$ is nondegenerate and that $t$
belongs to its interior relative to $J_{n+1}$. Then
$t_j\in F_{n+1,s}$ for all sufficiently large $j$.
Since $G_n|_{F_{n+1,s}}$ is continuous by its definition,
it follows that $\lim_{j\to\infty}G_n(t_j)=G_n(t)$.

It remains to consider the endpoints of $F_{n+1,s}$.
For terms $t_j\in F_{n+1,s}$, convergence follows from
continuity of the restriction. For terms outside $F_{n+1,s}$,
the points $r_{n+1}(t_j)$ approach $s$ from the side
corresponding to the endpoint $t$. Since $g_n$ has finitely
many nonconstant linear pieces, on each side sufficiently
close to $s$ its values approach $g_n(s)$ strictly from one
side. If these values are less than $g_n(s)$, then
\[
G_n(t_j)\le\ell_n(g_n(s)).
\]
Since every cluster point belongs to $F_{n,g_n(s)}$,
it follows that
\[
\lim_{j\to\infty}G_n(t_j)=\ell_n(g_n(s)).
\]
If the values are greater than $g_n(s)$, the same argument gives
\[
\lim_{j\to\infty}G_n(t_j)=u_n(g_n(s)).
\]

If $g_n$ is strictly increasing near $s$, the required left and right
endpoint values are $\ell_n(g_n(s))$ and $u_n(g_n(s))$,
respectively. These are the values given by $\psi_{n,s}=\id$.
If $g_n$ is strictly decreasing near $s$, the required values are
reversed, and are given by $\psi_{n,s}=R$.
If $s$ is a local maximum, both required endpoint values
are $\ell_n(g_n(s))$, and are given by $\psi_{n,s}=T$.
If $s$ is a local minimum, both required values are
$u_n(g_n(s))$, and are given by $\psi_{n,s}=R\circ T$.
The cases $s=0$ and $s=1$ follow from the same argument on
the one side that belongs to $I_{n+1}$. Therefore,
$G_n$ is continuous at every point of $J_{n+1}$.
By the definition of $G_n$, for each $s\in I_{n+1}$,
\begin{equation}\label{pa-eq:fibersonto}
G_n(F_{n+1,s})=F_{n,g_n(s)}.
\end{equation}
We prove that $G_n$ is surjective. Let $t\in J_n$.
Since $g_n$ is surjective, there is $s\in I_{n+1}$ such that
$g_n(s)=r_n(t)$. By \eqref{pa-eq:fibersonto},
\[
G_n(F_{n+1,s})=F_{n,g_n(s)},
\]
and $t\in F_{n,g_n(s)}$. Hence, $t\in G_n(J_{n+1})$.
Since $t\in J_n$ was arbitrary, $G_n$ is surjective.
The condition on the normalized restrictions follows directly
from the choice of $\psi_{n,s}$.
\end{proof}
We use the maps from Lemma~\ref{pa-lem:lifting} to define a map
between the corresponding inverse limits.

\begin{lemma}\label{pa-lem:limit-fibers}
Let $(G_n)$ be a sequence of maps as constructed in
Lemma~\ref{pa-lem:lifting}, let $X=\invlim(J_n,G_n)$, and, for each
positive integer $n$, let $q_n:X\to J_n$ be the $n$th coordinate
projection. Then $X$ is a chainable continuum and the map
$f:X\to Y$, defined by
\[
f((t_n))=(r_n(t_n))
\]
for each $(t_n)\in X$, is a continuous monotone surjection.
Moreover, for every $p\in Y$, every nondegenerate subcontinuum
of $f^{-1}(\{p\})$ contains an arc.
\end{lemma}
\begin{proof}
Since $X$ is an inverse limit of intervals with continuous
surjective bonding maps, it follows that $X$ is a chainable
continuum. By \eqref{pa-eq:commute}, the map
\[
f:X\to Y,\qquad f((t_n))=(r_n(t_n)),
\]
is well-defined. Since each $r_n$ is continuous, $f$ is
continuous. The commuting squares and the induced map are
shown in Figure~\ref{pa-fig:diagram}.

\begin{figure}[!ht]
\centering
\begin{tikzpicture}[font=\normalsize]
\node (X) at (0,1.7) {$X$};
\node (Jn) at (3.5,1.7) {$J_n$};
\node (Jnn) at (7.2,1.7) {$J_{n+1}$};
\node (P) at (0,0) {$Y$};
\node (In) at (3.5,0) {$I_n$};
\node (Inn) at (7.2,0) {$I_{n+1}$};
\draw[maparrow] (X)--node[left] {$f$} (P);
\draw[maparrow] (X)--node[above] {$q_n$} (Jn);
\draw[maparrow] (P)--node[below] {$\pi_n$} (In);
\draw[maparrow] (Jnn)--node[above] {$G_n$} (Jn);
\draw[maparrow] (Inn)--node[below] {$g_n$} (In);
\draw[maparrow] (Jn)--node[right] {$r_n$} (In);
\draw[maparrow] (Jnn)--node[right] {$r_{n+1}$} (Inn);
\node at (9,1.7) {$\cdots$};
\node at (9,0) {$\cdots$};
\draw[maparrow] (8.6,1.7)--(Jnn);
\draw[maparrow] (8.6,0)--(Inn);
\end{tikzpicture}
\caption{The maps between the inverse sequences and their inverse limits.}
\label{pa-fig:diagram}
\end{figure}
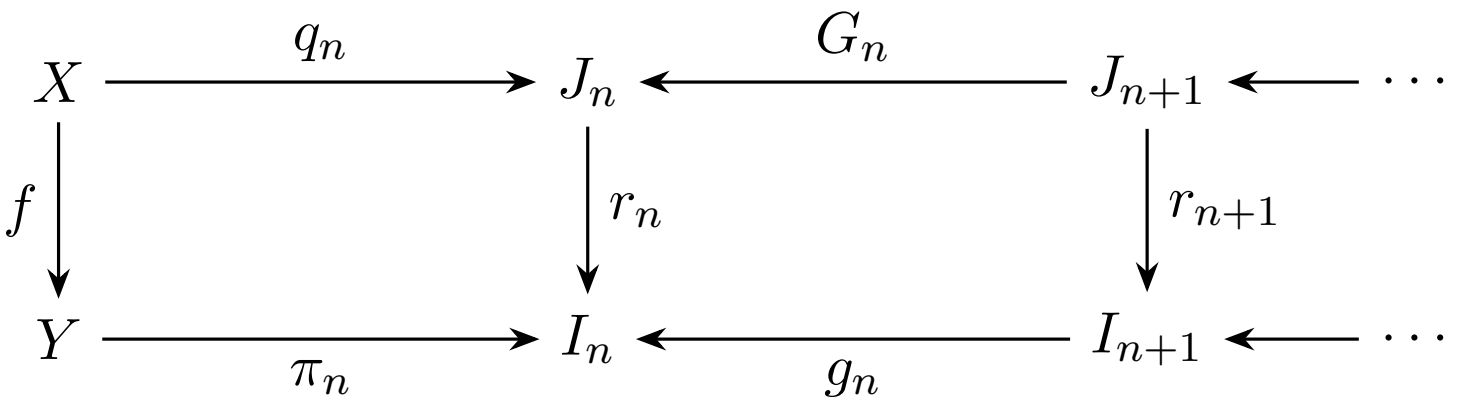

Let $p=(t_n)\in Y$. We prove that $f^{-1}(\{p\})$ is a
non-empty continuum. For each positive integer $n$, we have
$g_n(t_{n+1})=t_n$. Therefore, by Lemma~\ref{pa-lem:lifting},
\[
G_n(F_{n+1,t_{n+1}})=F_{n,t_n},
\]
so the restrictions
$G_n|_{F_{n+1,t_{n+1}}}:F_{n+1,t_{n+1}}\to F_{n,t_n}$
form an inverse sequence.

A point $(s_n)\in X$ belongs to $f^{-1}(\{p\})$ if and only if
$r_n(s_n)=t_n$ for each positive integer $n$. Since
$r_n^{-1}(\{t_n\})=F_{n,t_n}$, this is equivalent to
$s_n\in F_{n,t_n}$ for each positive integer $n$.
Together with the equations $G_n(s_{n+1})=s_n$ defining $X$,
these are precisely the conditions for membership in
$\invlim(F_{n,t_n},G_n|_{F_{n+1,t_{n+1}}})$. Hence,
\begin{equation}\label{pa-eq:fiberlimit}
f^{-1}(\{p\})=
\invlim\bigl(F_{n,t_n},G_n|_{F_{n+1,t_{n+1}}}\bigr).
\end{equation}
Every restricted bonding map in \eqref{pa-eq:fiberlimit}
is a continuous surjection between non-empty compact intervals.
Therefore, $f^{-1}(\{p\})$ is a non-empty continuum.
Since $p\in Y$ was arbitrary, $f$ is surjective and monotone.

We next prove that for each $p$, each nondegenerate subcontinuum of $f^{-1}(\{p\})$ contains an arc.
By Lemma~\ref{pa-lem:dense-sets},
$g_n^{-1}(D_n)=D_{n+1}$ for each positive integer $n$.
Since $g_n(t_{n+1})=t_n$, it follows that
\[
t_{n+1}\in D_{n+1}\quad\Longleftrightarrow\quad t_n\in D_n
\]
for each positive integer $n$. Therefore, either
$t_n\notin D_n$ for every positive integer $n$, or
$t_n\in D_n$ for every positive integer $n$.
In the first case, every interval $F_{n,t_n}$ is a singleton.
By \eqref{pa-eq:fiberlimit}, the fiber $f^{-1}(\{p\})$ is also
a singleton, so it has no nondegenerate subcontinua.
Suppose that the second case holds. For each positive integer
$n$, let
\[
E_n=F_{n,t_n},\qquad
b_n=a_{n,t_n},\qquad
h_n=G_n|_{F_{n+1,t_{n+1}}}.
\]
Then each $E_n$ is a nondegenerate compact interval, and $b_n$
is the increasing affine homeomorphism from $[0,1]$ onto $E_n$.
Since $g_n(t_{n+1})=t_n$, Lemma~\ref{pa-lem:lifting} gives
\[
b_n^{-1}\circ h_n\circ b_{n+1}
=\psi_{n,t_{n+1}}\in\{\id,R,T,R\circ T\}.
\]
Therefore, by Lemma~\ref{pa-lem:arcs} and \eqref{pa-eq:fiberlimit},
every nondegenerate subcontinuum of $f^{-1}(\{p\})$
contains an arc. This proves the assertion in both cases.
\end{proof}

The following lemma is the last step needed to prove our main result.

\begin{lemma}\label{pa-lem:arc-extraction}
Let $X$ and $f:X\to Y$ be as in Lemma~\ref{pa-lem:limit-fibers}.
Then every nondegenerate subcontinuum of $X$ contains an arc.
\end{lemma}
\begin{proof}
Let $H$ be a nondegenerate subcontinuum of $X$. We prove that
$H$ contains an arc. If $f(H)$ is degenerate, then
$H\subseteq f^{-1}(\{p\})$ for some $p\in Y$. In this case,
the assertion follows from Lemma~\ref{pa-lem:limit-fibers}.

Suppose that $f(H)$ is nondegenerate. Then there is a positive
integer $n$ such that $\pi_n(f(H))$ is nondegenerate. Let $n$ be
such an integer. Since $D_n$ is dense in $I_n$, there is
\[
d\in D_n\cap\Int\bigl(\pi_n(f(H))\bigr).
\]
Let $d$ be such a point. Note that $q_n(H)$ is a compact interval
and that
\[
r_n(q_n(H))=\pi_n(f(H)).
\]
Since $d$ belongs to the interior of this image, the interval
$q_n(H)$ contains points on both sides of $F_{n,d}$. Therefore,
$F_{n,d}\subseteq q_n(H)$. Also, $F_{n,d}$ is nondegenerate
because $d\in D_n$.

Consider the continuous surjection $q_n|_H:H\to q_n(H)$.
By Lemma~\ref{pa-lem:crossing}, there is a subcontinuum $L$ of $H$
such that $q_n(L)=F_{n,d}$. Let $L$ be such a continuum.
Since $F_{n,d}$ is nondegenerate, $L$ is nondegenerate. Moreover,
\[
\pi_n(f(L))=r_n(q_n(L))=r_n(F_{n,d})=\{d\}.
\]
It follows that $f(L)\subseteq\pi_n^{-1}(\{d\})$.
By Lemma~\ref{pa-lem:coordinate-fibers}, $\pi_n^{-1}(\{d\})$ is
zero-dimensional. Since $f(L)$ is connected, it is a singleton.
Let $p\in Y$ be such that $f(L)=\{p\}$. Then
$L\subseteq f^{-1}(\{p\})$. By Lemma~\ref{pa-lem:limit-fibers},
$L$ contains an arc. Since $L\subseteq H$, it follows that $H$
contains an arc, as required.
\end{proof}

\subsubsection{A universal domain containing no pseudoarc}\label{pa-sec:main}

We use the preceding lemmas
to construct a chainable reciprocating domain containing no pseudoarc.

\begin{theorem}\label{pa-thm:counterexample}
Let $P$ be a pseudoarc. There exist a chainable continuum $X$ and a
continuous monotone surjection $f:X\to P$ such that every nondegenerate
subcontinuum of $X$ contains an arc. In particular, $X$ contains no
subcontinuum homeomorphic to $P$.
\end{theorem}

\begin{proof}
Let $P$ be a pseudoarc. By the inverse-limit representation
recalled in \cite[Section 1.3]{AAH}, $P$ admits an inverse-limit
representation with nowhere constant piecewise linear bonding
maps. So, suppose that $P$ is homeomorphic to such an inverse limit $Y=\invlim(I_n,g_n)$ and let $\eta:Y\to P$ be a homeomorphism.

By Lemma~\ref{pa-lem:dense-sets}, let $(D_n)$ be a sequence of countable
dense sets with the stated properties. By Lemma~\ref{pa-lem:expansion},
let $(r_n)$ be the corresponding sequence of continuous monotone
surjections, and, by Lemma~\ref{pa-lem:lifting}, let $(G_n)$ be the
corresponding sequence of lifted bonding maps. Let
$X=\invlim(J_n,G_n)$. By Lemma~\ref{pa-lem:limit-fibers}, $X$ is a
chainable continuum and there is a continuous monotone surjection $f_0:X\to Y$.
Then $f=\eta\circ f_0:X\to P$ is a continuous monotone surjection.
By Lemma~\ref{pa-lem:arc-extraction}, every nondegenerate subcontinuum
of $X$ contains an arc. Since a pseudoarc contains no arc, $X$
contains no subcontinuum homeomorphic to $P$.
\end{proof}

\begin{corollary}\label{cor:domain-no-pseudoarc}
There is a reciprocating and universal domain $X$ in
$\mathcal{CC}$ such that every nondegenerate subcontinuum of $X$
contains an arc. In particular, $X$ contains no nondegenerate
hereditarily indecomposable subcontinuum. The same $X$ is a
reciprocating and universal domain in $\mathcal{WC}$.
\end{corollary}
\begin{proof}
Choose $X$ as in Theorem~\ref{pa-thm:counterexample}. Since $X$
maps onto a pseudoarc, Theorems~\ref{DDJI} and \ref{WC} imply
that it is reciprocating and universal in both classes. Every
nondegenerate hereditarily indecomposable continuum contains no
arc, so the assertion about subcontinua follows.
\end{proof}

Although the pseudoarc is universal in the class of chainable continua, it is not universal in the broader class of tree-like continua. This is because every continuous image of the pseudoarc is weakly chainable, yet there exist tree-like continua that are not weakly chainable, for an example see \cite{moebes}. Consequently, such continua cannot be continuous images of the pseudoarc. Therefore, the pseudoarc is not a universal domain for the class of tree-like continua. We continue studying tree-like continua in Subsection \ref{treelike}.
\subsection{Tree-like continua}\label{treelike}

In this subsection, we prove that there are no reciprocating domains in
the classes of tree-like continua and planar tree-like continua. We first
show that tree-likeness is preserved under wedges and then combine this
with known nonexistence results for universal tree-like continua. 

\begin{definition}
	A {\em graph\/} is a singleton or a continuum that can be written as the union of finitely many arcs (continua that are homeomorphic to $[0,1]$), any two of which are either disjoint or intersect in one or both of their end points. A graph is {\em acyclic\/} if it does not contain a simple closed curve. An acyclic graph is called a {\em tree\/}. 
\end{definition} 

\begin{definition}
	Let $X$ be a continuum. We say that $X$ is a \emph{tree-like continuum} if there is an inverse sequence $(T_{\ell},f_{\ell})$, where for each positive integer $\ell$, $T_{\ell}$ is a tree, such that $X$ is homeomorphic to the inverse limit $\invlim(T_\ell,f_\ell)$.
	\end{definition}
	\begin{definition}
	We denote by $\mathcal{TLC}$ the class of all tree-like continua.
\end{definition}
Moreover, tree-likeness is preserved under inverse limits of tree-like continua, as shown by the following result of Marde\v si\' c.

\begin{theorem}\label{tl}
	Let $(X_{\ell},f_{\ell})$ be an inverse sequence, where for each positive integer $\ell$, $X_{\ell}$ is a tree-like continuum. Then the inverse limit $\invlim(X_\ell,f_\ell)$ is a tree-like continuum.
\end{theorem}

\begin{proof}
	The proof can be found in  \cite{2}.
\end{proof}
\begin{observation}\label{wedgetree}
	Note that any wedge of two trees is a tree. 
\end{observation}
Next, we show that wedges preserve tree-likeness. This is a
special case of a theorem of Cook \cite[Theorem~2]{Cook1970},
which states that the union of two tree-like continua with
nonempty connected intersection is tree-like; see also
\cite[Lemma~6(1)]{HagopianMarsh2015}.
We include a direct proof using inverse limits of trees.

\begin{theorem}\label{wedgetreelike}
Let $X=\invlim(T_n,f_n)$ and $Y=\invlim(S_n,g_n)$ be inverse
limits of trees, and let $t=(t_n)\in X$ and $s=(s_n)\in Y$.
Then $W=(X\sqcup Y)/_{\{t,s\}}$ is a tree-like continuum.
\end{theorem}

\begin{proof}
For each positive integer $n$, let
\[
W_n=(T_n\sqcup S_n)/_{\{t_n,s_n\}},
\]
and define $h_n:W_{n+1}\to W_n$ by
\[
h_n([z])=
\begin{cases}
[f_n(z)],&z\in T_{n+1},\\
[g_n(z)],&z\in S_{n+1}.
\end{cases}
\]
The definitions agree at the identified points, since
$f_n(t_{n+1})=t_n$ and $g_n(s_{n+1})=s_n$.
The quotient property and the pasting lemma
\cite[Theorem 18.3]{Munkres} imply continuity of $h_n$.
By Observation~\ref{wedgetree}, each $W_n$ is a tree.
The coordinate quotient maps induce a continuous map
\[
\Phi:W\longrightarrow\invlim(W_n,h_n),\qquad
\Phi([(z_n)])=([z_n]).
\]
To see that $\Phi$ is injective, suppose that two points
have the same image. If both belong to $X$ or both
belong to $Y$, their coordinates agree, so they are equal.
If one belongs to $X$ and the other to $Y$, equality
of their images forces every coordinate to be an
identified point. Hence the two points are $t$ and $s$,
which represent the same point of $W$.

To prove surjectivity, let $(w_n)\in\invlim(W_n,h_n)$.
If every $w_n$ is the identified point of $W_n$,
then $(w_n)=\Phi([t])$.
Otherwise, choose $n$ such that $w_n$ is not the
identified point. Then $w_n$ belongs to exactly one
of $T_n$ and $S_n$, and the bonding maps force every
subsequent coordinate to belong to the same tree
sequence. The preceding coordinates are determined
by the bonding maps. Thus $(w_n)$ is the image under
$\Phi$ of a point of $X$ or $Y$.

Consequently, $\Phi$ is a continuous bijection from
a compact space onto a Hausdorff space and therefore
a homeomorphism. Since each $W_n$ is a tree, $W$ is
tree-like.
\end{proof}

Next, we prove that there are no universal domains in $\mathcal{TLC}$. We start with the following two results by Ingram and Russo.
\begin{theorem}[Ingram]\label{ingi}
	For any $X\in \mathcal{TLC}$, there is $Y\in \mathcal{TLC}$ such that there are no continuous surjections from $X$ to $Y$.
\end{theorem}
\begin{proof}
	In \cite{ingram}, a class $\mathcal G$ of tree-like continua such that for any continuum $M$, there is a member of $\mathcal G$ which is not a continuous image of $M$ is constructed. Let $X\in \mathcal{TLC}$. Then there is $Z\in \mathcal G$ such that there are no continuous surjections from $X$ to $Z$. Let $Y\in \mathcal{TLC}$ be such that $Z$ is homeomorphic to $Y$.  Then there are no continuous surjections from $X$ to $Y$. 
\end{proof}

\begin{corollary}\label{taonki}
	For each $X\in \mathcal{TLC}$, $X$ is not a universal domain in $\mathcal{TLC}$.
\end{corollary}
\begin{proof}
	Let $X\in \mathcal{TLC}$ and suppose that $X$ is a universal domain in $\mathcal{TLC}$. Then for each $Y\in \mathcal{TLC}$, there is a continuous surjection from $X$ to $Y$.  This contradicts Theorem \ref{ingi}. 
\end{proof}
The following theorem is the main result of this section.
\begin{theorem}
	For each $X\in \mathcal{TLC}$, $X$ is not a reciprocating domain in $\mathcal{TLC}$.
\end{theorem}
\begin{proof}
	Let $X\in \mathcal{TLC}$ and suppose that $X$ is a reciprocating domain in $\mathcal{TLC}$. Since $\mathcal{TLC}$ is closed for wedges by Theorem \ref{wedgetreelike}, it follows from Corollary \ref{wedge1} that $X$ is a universal domain, which contradicts Corollary \ref{taonki}. 
\end{proof}

We continue with planar tree-like continua.
\begin{definition}
	We use $\mathcal{PTLC}$ to denote the class of all planar tree-like continua. 
\end{definition}

\begin{theorem}[Russo]\label{rus}
	For each $X\in \mathcal{PTLC}$, $X$ is not a universal domain in $\mathcal{PTLC}$.
\end{theorem}
\begin{proof}
	Theorem \ref{rus} is (in its equivalent form) one of the main results in \cite{russo}.
	\end{proof}
The following theorem is also one of the main results of this section.
\begin{theorem}
	For each $X\in \mathcal{PTLC}$, $X$ is not a reciprocating domain in $\mathcal{PTLC}$.
\end{theorem}
\begin{proof}
	Let $X\in \mathcal{PTLC}$ and suppose that $X$ is a reciprocating domain in $\mathcal{PTLC}$. Note that for any two planar tree-like continua $X$ and $Y$, there is a planar tree-like continuum that is a wedge of $X$ and $Y$. Therefore, it follows from Corollary \ref{wedge1} that $X$ is a universal domain in $\mathcal{PTLC}$, which contradicts Theorem \ref{rus}. 
\end{proof}

\subsection{Solenoids}

In this subsection, we study continuous surjections between solenoids and
characterize reciprocating domains in the class of solenoids. In
particular, we relate reciprocating solenoids to universal solenoids.
\begin{definition}
We use $S^1$ to denote the unit circle in the complex plane; $S^1=\{z\in \mathbb C \mid |z|=1\}$. 
For each prime $p$, we denote by $f_p:S^1\to S^1$ the map that is defined by $f_p(z)=z^p$ for each $z\in S^1$. 
\end{definition}

\begin{definition}\label{altu}
Let $X$ be a continuum. We say that $X$ is a \emph{solenoid} if it is homeomorphic to the inverse limit
\[
\varprojlim \bigl(S^1,f_{p_i}\bigr)_{i=1}^{\infty}
\]
for some sequence $(p_i)$ of primes. In this case we say that $X$ is a $(p_i)$-adic solenoid. A solenoid is called a \emph{universal solenoid} if every prime occurs infinitely many times in the sequence $(p_i)$.
\end{definition}
\begin{definition}
 We denote by $\mathcal S$ the class of all solenoids in $\mathcal C$.  
\end{definition}
\begin{definition}
	Let $\mathbb N$ be the set of positive integers  and let $P = (p_1, p_2, p_3,\ldots )$ and
$Q = (q_1, q_2, q_3,\ldots)$ be any sequences of primes.  Then we define $P \unlhd Q$ as follows:  $P \unlhd Q$ if and only if there is a set $N\subseteq \mathbb N$ such that $\mathbb N\setminus N$ is finite and there is an injective function $\alpha:N\rightarrow \mathbb N$ such that for each $n\in N$,
$$
p_n=q_{\alpha(n)}.
$$
\end{definition}
The following three theorems are well-known results; the reader can find them in Prajs's paper \cite{prajs}. They are stated (and proved) by Cook and McCord in equivalent forms in \cite{Cook1967} and \cite{McCord1965}. We use them to prove Theorem \ref{pr}.
\begin{theorem}\label{prasje1}
	Let $P=(p_1, p_2, p_3,\ldots )$ and $Q = (q_1, q_2, q_3,\ldots)$ be any sequences of primes and let $X,Y\in \mathcal S$ be such that $X$ is homeomorphic to $\varprojlim \bigl(S^1,f_{p_i}\bigr)_{i=1}^{\infty}$ and $Y$ is homeomorphic to $\varprojlim \bigl(S^1,f_{q_i}\bigr)_{i=1}^{\infty}$. The following statements are equivalent.
	\begin{enumerate}
		\item There is a continuous surjection from $Y$ to $X$.
		\item $P\unlhd Q$.
	\end{enumerate}
\end{theorem}
\begin{proof}
	See \cite[Theorem 4.4, page 348]{prajs}.
\end{proof}
\begin{theorem}\label{prasje2}
	Let $P=(p_1, p_2, p_3,\ldots )$ and $Q = (q_1, q_2, q_3,\ldots)$ be any sequences of primes and let $X,Y\in \mathcal S$ be such that $X$ is homeomorphic to $\varprojlim \bigl(S^1,f_{p_i}\bigr)_{i=1}^{\infty}$ and $Y$ is homeomorphic to $\varprojlim \bigl(S^1,f_{q_i}\bigr)_{i=1}^{\infty}$. The following statements are equivalent.
	\begin{enumerate}
		\item $Y$ is homeomorphic to $X$.
		\item $Q\unlhd P$ and $P\unlhd Q$.
	\end{enumerate}
\end{theorem}
\begin{proof}
	See \cite[Theorem 4.5, page 348]{prajs}.
\end{proof}

\begin{corollary}\label{prasje}
	Let $X,Y\in \mathcal S$. The following statements are equivalent.
	\begin{enumerate}
		\item $X$ is homeomorphic to  $Y$.
		\item There is a continuous surjection from $X$ to $Y$ and there is a continuous surjection from $Y$ to $X$.
	\end{enumerate}  
\end{corollary}
\begin{proof}
	See \cite[Theorem 4.6, page 349]{prajs}.
\end{proof}
The next theorem is a well known result, but we provide a brief proof for completeness.
\begin{theorem}\label{universal-unique}
Any two universal solenoids are homeomorphic.
\end{theorem}

\begin{proof}
Let $X$ and $Y$ be any universal solenoids, with defining prime sequences
$P=(p_1,p_2,p_3,\ldots)$ and $Q=(q_1,q_2,q_3,\ldots)$ respectively. 
By universality, every prime appears infinitely often in $P$ and also in $Q$.
First we show $Q\unlhd P$. Since every prime that occurs in $Q$ occurs infinitely many times in $P$,
we can choose, for each $n$, a distinct index $\alpha(n)$ with $q_n=p_{\alpha(n)}$ in such a way that this defines an injection $\alpha:\mathbb N\to\mathbb N$ such that $q_n=p_{\alpha(n)}$ for all $n\in \mathbb N$. Hence $Q\unlhd P$.
By symmetry, the same reasoning gives $P\unlhd Q$.
Now apply Theorem~\ref{prasje1} and Corollary~\ref{prasje}: 
from $Q\unlhd P$ and $P\unlhd Q$ we obtain continuous surjections from $Y$ to $X$ and from $X$ to $Y$;
then Corollary~\ref{prasje} yields that $X$ and $Y$ are homeomorphic.
\end{proof}
We use the above theorem to prove the following main result of this section, where a characterization of reciprocating domains in $\mathcal S$ is obtained.

\begin{theorem}\label{pr}
For each $X\in \mathcal S$, the following statements are equivalent:
\begin{enumerate}
\item\label{S1} $X$ is a reciprocating domain in $\mathcal S$.
\item\label{S3} $X$ is the universal solenoid.
\item\label{S2} $X$ is a universal domain in $\mathcal S$.
\end{enumerate}
\end{theorem}

\begin{proof}
Let $U$ be a solenoid whose defining prime sequence contains each
prime infinitely often. For any solenoid $Y$ with defining sequence
$Q$, there is an injection matching every term of $Q$ to an equal
term in the sequence of $U$. By Theorem~\ref{prasje1}, $U$ maps
continuously onto $Y$. Thus $U$ is a universal domain in $\mathcal S$.
Theorem~\ref{zelosplosno} shows that statements \ref{S1} and
\ref{S2} are equivalent.

If $X$ is universal, it maps onto the $(p,p,p,\ldots)$-adic
solenoid for every prime $p$. By Theorem~\ref{prasje1}, each
prime therefore occurs infinitely often in the defining sequence of
$X$. Hence $X$ is a universal solenoid. Conversely, the argument
used for $U$ applies to every universal solenoid. This proves the
equivalence with statement \ref{S3}.
\end{proof}
\subsection{Circle-like continua}

In this subsection, we characterize reciprocating domains among
circle-like continua. We use a universal circle-like continuum to show
that reciprocating and universal domains coincide in this class.

\begin{definition}
	Let $X$ be a continuum. We say that $X$ is a \emph{circle-like continuum} if there is an inverse sequence $(C_{\ell},f_{\ell})$, where for each positive integer $\ell$, $C_{\ell}$ is a circle, such that $X$ is homeomorphic to the inverse limit $\invlim(C_\ell,f_\ell)$.
\end{definition}

\begin{observation}
	Note that solenoids are examples of circle-like continua. 
\end{observation}

\begin{definition}
	We denote by $\mathcal{CLC}$ the class of all circle-like continua.
\end{definition}

\begin{theorem}\label{main:P1}
There exists $X \in \mathcal{CLC}$ such that $X$ is a universal as well as a reciprocating domain in $\mathcal{CLC}$.
\end{theorem}

\begin{proof}
In \cite[Corollary 5.30]{bartos}, Barto\v s and Kubi\' s show that there is a circle-like continuum, the so-called universal pseudo-solenoid $\mathbb P_{\prod}$, which is a universal domain in $\mathcal{CLC}$. It follows from Observation~\ref{karin} that the universal pseudo-solenoid $\mathbb P_{\prod}$ is a reciprocating domain in $\mathcal{CLC}$.
\end{proof}

\begin{theorem}\label{CLC}
	Let $X\in \mathcal{CLC}$ and let $\mathbb P_{\prod}$ be the universal pseudo-solenoid from \cite[Corollary 5.30]{bartos}. The following statements are equivalent.
	\begin{enumerate}
		\item\label{juc11d} $X$ is a reciprocating domain in $\mathcal{CLC}$.
		\item\label{juc21d} There is a continuous surjection $f:X\to \mathbb P_{\prod}$.
		\item\label{juc31d} $X$ is a universal domain in $\mathcal{CLC}$.
	\end{enumerate}
\end{theorem}

\begin{proof}
	The theorem is a corollary to Theorem \ref{zelosplosno}.
\end{proof}

\subsection{Uniformly path connected continua}

In this subsection, we use the Cantor fan as a universal domain for
uniformly path connected continua. We characterize reciprocating domains
in this class as precisely those continua that admit a continuous
surjection onto the Cantor fan.
\begin{definition}
	Let $X$ be a continuum. We say that $X$ is \emph{uniformly path connected} if there exists a family
$\mathcal P$ of paths $p:[0,1]\to X$ in $X$ such that
\begin{enumerate}
	\item for all $x,y\in X$ there is a path in $\mathcal P$ from $x$ to $y$, and
	\item for any $\varepsilon>0$ there is a positive integer $n$ such that for each $p\in \mathcal P$, there are $t_0$, $t_1$, $t_2$, $\ldots$, $t_n\in [0,1]$ such that 
	\begin{enumerate}
	\item $t_0=0$ and $t_n=1$,
	\item for each $i\in\{0,1,2,\ldots, n-1\}$, $t_{i}<t_{i+1}$, and
	\item for each $i\in\{0,1,2,\ldots, n-1\}$, $\diam\Big(p\big([t_i,t_{i+1}]\big)\Big)\leq \varepsilon$.
	\end{enumerate} 
\end{enumerate}
\end{definition}

\begin{definition}
Let $X$ be a continuum. We say that $X$ is \emph{a Cantor fan} if $X$ is homeomorphic to the continuum $\bigcup_{c\in C}A_c$, where $C\subseteq [0,1]$ is the Cantor middle third set and for each $c\in C$, $A_c$ is the  {convex} segment in the plane from $(\frac{1}{2},0)$ to $(c,1)$.  See Figure \ref{fig:cantor-fan}.
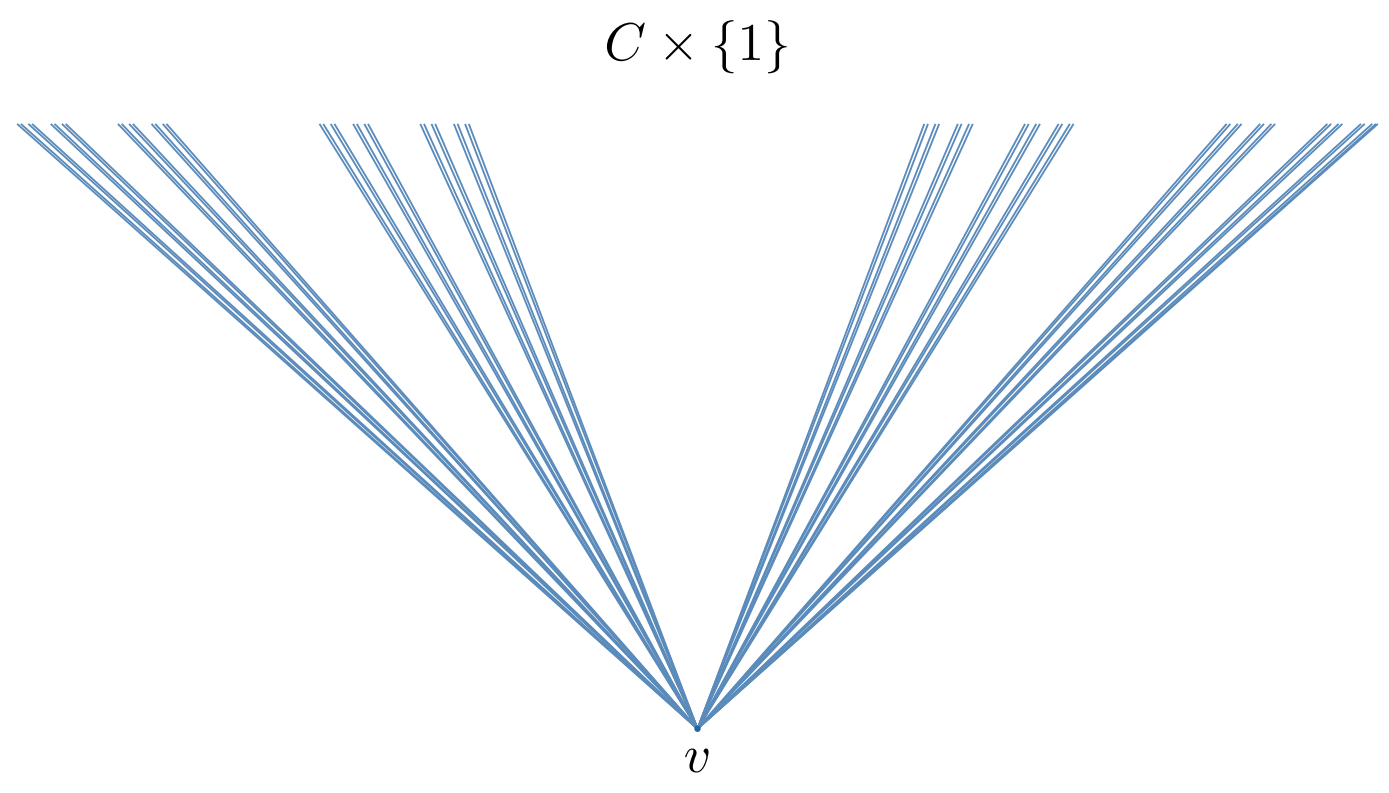
\begin{figure}[h!]
	\centering
		\begin{tikzpicture}[x=9cm,y=4cm,font=\small]
\coordinate (v) at (.5,0);
\foreach \a in {0,2}
\foreach \b in {0,2}
\foreach \c in {0,2}
\foreach \d in {0,2}
\foreach \e in {0,2}
\foreach \f in {0,2}{
  \pgfmathsetmacro{\xx}{\a/3+\b/9+\c/27+\d/81+\e/243+\f/729}
  \draw[fiberblue!75,line width=.35pt] (v)--(\xx,1);
}
\draw[fiberblue!75,line width=.35pt] (v)--(1,1);
\fill[fiberblue] (v) circle (.6pt);
\node[below] at (v) {$v$};
\node[above] at (.5,1.05) {$C\times\{1\}$};
\end{tikzpicture}
	\caption{A drawing of the Cantor fan. Only finitely many legs are shown.}
	\label{fig:cantor-fan}
\end{figure}  
\end{definition}
The following is the main result from \cite{kuperberg}.
\begin{theorem}[Kuperberg]\label{kup}
	Let $X$ be a continuum. The following statements are equivalent.
	\begin{enumerate}
		\item $X$ is uniformly path connected.
		\item There is a continuous surjection from the Cantor fan to $X$.
	\end{enumerate}
\end{theorem}
\begin{proof}
	Theorem follows directly from \cite[Theorem 3.5, page 322]{kuperberg}.
\end{proof}
\begin{definition}
	We use $\mathcal{UPC}$ to denote the class of all uniformly path connected continua. 
\end{definition}
\begin{observation}
Let $C$ be a Cantor fan. Since the identity map from $C$ to $C$ is a continuous surjection, it follows that $C\in \mathcal{UPC}$.
\end{observation}
\begin{corollary}\label{kupi}
	Let $C$ be a Cantor fan. Then $C$ is a universal domain in $\mathcal{UPC}$.
\end{corollary}
\begin{proof}
	The corollary follows directly from Theorem \ref{kup}.
\end{proof}
\begin{theorem}\label{UPC}
	Let $X\in \mathcal{UPC}$ and let $C$ be a Cantor fan. The following statements are equivalent.
	\begin{enumerate}
		\item\label{juc11} $X$ is a reciprocating domain in $\mathcal{UPC}$.
		\item\label{juc21} There is a continuous surjection $f:X\to C$.
		\item\label{juc31} $X$ is a universal domain in $\mathcal{UPC}$.
	\end{enumerate}
\end{theorem}

\begin{proof}
	The theorem is a corollary to Theorem \ref{zelosplosno}.
\end{proof}

\subsection{Smooth dendroids}

In this subsection, we use the Cantor fan as a universal domain for smooth
dendroids. Consequently, we characterize reciprocating domains among
smooth dendroids by the existence of a continuous surjection onto the
Cantor fan.
 \begin{definition}
Let $X$ be a continuum. 
\begin{enumerate}
\item The continuum $X$ is \emph{unicoherent} if for any subcontinua $A$ and $B$ of $X$ such that $X=A\cup B$,  the compactum $A\cap B$ is connected. 
\item The continuum $X$ is \emph{hereditarily unicoherent } provided that each of its subcontinua is unicoherent.
\item The continuum $X$ is a \emph{dendroid} if it is an arcwise connected, hereditarily unicoherent continuum.
\item Let $X$ be a dendroid. For all points $x$ and $y$ in $X$, we define  \emph{$[x,y]$} to be the arc in $X$ with end points $x$ and $y$, if $x\neq y$. If $x=y$, then we define $[x,y]=\{x\}$.
\item Let $X$ be a dendroid. We say that the dendroid $X$ is \emph{smooth} if there is $p\in X$ such that for any $x\in X$ and for any sequence $(x_n)$ of points in $X$,
$$
\lim_{n\to \infty}x_n=x \Longrightarrow \lim_{n\to \infty}[p,x_n]=[p,x].
$$ 
\end{enumerate}
\end{definition}
\begin{definition}
	We use $\mathcal{SD}$ to denote the class of all smooth dendroids.
\end{definition}
\begin{observation}
	Let $C$ be a Cantor fan. Note that $C\in \mathcal{SD}$.
\end{observation}
\begin{theorem}[Charatonik, Eberhart]\label{CE}
	Let $C$ be a Cantor  fan and let $X$ be a continuum. If $X\in \mathcal{SD}$, then there is a continuous surjection from $C$ to $X$.
\end{theorem}
\begin{proof}
	The theorem follows from \cite[Corollary 11, page 311]{charatonike}.
\end{proof}
\begin{corollary}\label{lidi}
	Let $C$ be a Cantor fan. Then $C$ is a universal domain in $\mathcal{SD}$. 
\end{corollary}
\begin{proof}
	The corollary follows directly from Theorem \ref{CE}.
\end{proof}
\begin{theorem}\label{SD}
	Let $X\in \mathcal{SD}$  and let $C$ be a Cantor fan. The following statements are equivalent.
	\begin{enumerate}
		\item\label{juc111} $X$ is a reciprocating domain in $\mathcal{SD}$.
		\item\label{juc211} There is a continuous surjection $f:X\to C$.
		\item\label{juc311} $X$ is a universal domain in $\mathcal{SD}$.
	\end{enumerate}
\end{theorem}

\begin{proof}
	The theorem is a corollary to Theorem \ref{zelosplosno}.
\end{proof}
\subsection{Smooth fans}

In this subsection, we study reciprocating domains in the class of smooth
fans. We show that both the Cantor fan and the Lelek fan are universal
domains and obtain several equivalent characterizations of reciprocating
smooth fans in terms of continuous surjections, contained subfans, and
the number of legs.
 \begin{definition}
Let $X$ be a continuum. 
\begin{enumerate} 
\item Let $X$ be an arc. A point $x\in X$ is called and \emph{end point of $X$} if  there is a homeomorphism $\varphi:[0,1]\rightarrow X$ such that $\varphi(0)=x$.
\item Let $X$ be a dendroid.  A point $x\in X$ is called an \emph{end point of $X$} if for  every arc $A$ in $X$ that contains $x$, $x$ is an end point of $A$.  The set of all end points of $X$ is denoted by $E(X)$. 
\item The continuum $X$ is a \emph{simple triod} if it is homeomorphic to $([-1,1]\times \{0\})\cup (\{0\}\times [0,1])$.
\item Let $X$ be a simple triod. A point $x\in X$ is called the \emph{top point} or, briefly, \emph{the top of $X$} if  there is a homeomorphism $\varphi:([-1,1]\times \{0\})\cup (\{0\}\times [0,1])\rightarrow X$ such that $\varphi(0,0)=x$.
\item Let $X$ be a dendroid.  A point $x\in X$ is called a \emph{ramification point of $X$}, if there is a simple triod $T$ in $X$ with top $x$.  The set of all ramification points of $X$ is denoted by $R(X)$. 
\item The continuum $X$ is a \emph{fan} if it is a dendroid with at most one ramification point $v$, which is called \emph{the top of the fan $X$} (if it exists).
\item Let $X$ be a fan with top $v$. For each end point $e\in E(X)$ of the fan $X$, we call the arc in $X$ from $v$ to $e$, a \emph{leg of $X$}. The set of all legs of $X$ is denoted by $\mathcal L(X)$.
\item Let $X$ be a fan with top $v$. We say that the fan $X$ is \emph{smooth} if for any $x\in X$ and for any sequence $(x_n)$ of points in $X$,
$$
\lim_{n\to \infty}x_n=x \Longrightarrow \lim_{n\to \infty}[v,x_n]=[v,x].
$$ 
\item Let $X$ be a fan.    We say that $X$ is a \emph{Lelek fan} if it is nondegenerate and smooth and $\Cl(E(X))=X$. See Figure \ref{figure2}.
\begin{figure}[h!]
	\centering
		\begin{tikzpicture}[x=9cm,y=4cm,font=\small]
\coordinate (v) at (.5,0);
\foreach \k in {0,...,127}{
  \pgfmathsetmacro{\xx}{2*mod(\k,2)/3
    +2*mod(floor(\k/2),2)/9
    +2*mod(floor(\k/4),2)/27
    +2*mod(floor(\k/8),2)/81
    +2*mod(floor(\k/16),2)/243
    +2*mod(floor(\k/32),2)/729
    +2*mod(floor(\k/64),2)/2187}
  \pgfmathsetmacro{\hh}{.06+.94*mod(53*\k+19,128)/127}
  \pgfmathsetmacro{\xt}{.5+\hh*(\xx-.5)}
  \draw[fiberblue!65,line width=.3pt] (v)--(\xt,\hh);
  \fill[branchgreen] (\xt,\hh) circle (.35pt);
}
\fill[fiberblue] (v) circle (.6pt);
\node[below] at (v) {$v$};
\node[above,align=center] at (.5,1.08)
  {Endpoints at different heights};
\end{tikzpicture}
	\caption{A drawing of the Lelek fan. The finite collection of legs illustrates endpoints at varying heights; the endpoint set of the Lelek fan itself is dense.}
	\label{figure2}
\end{figure}
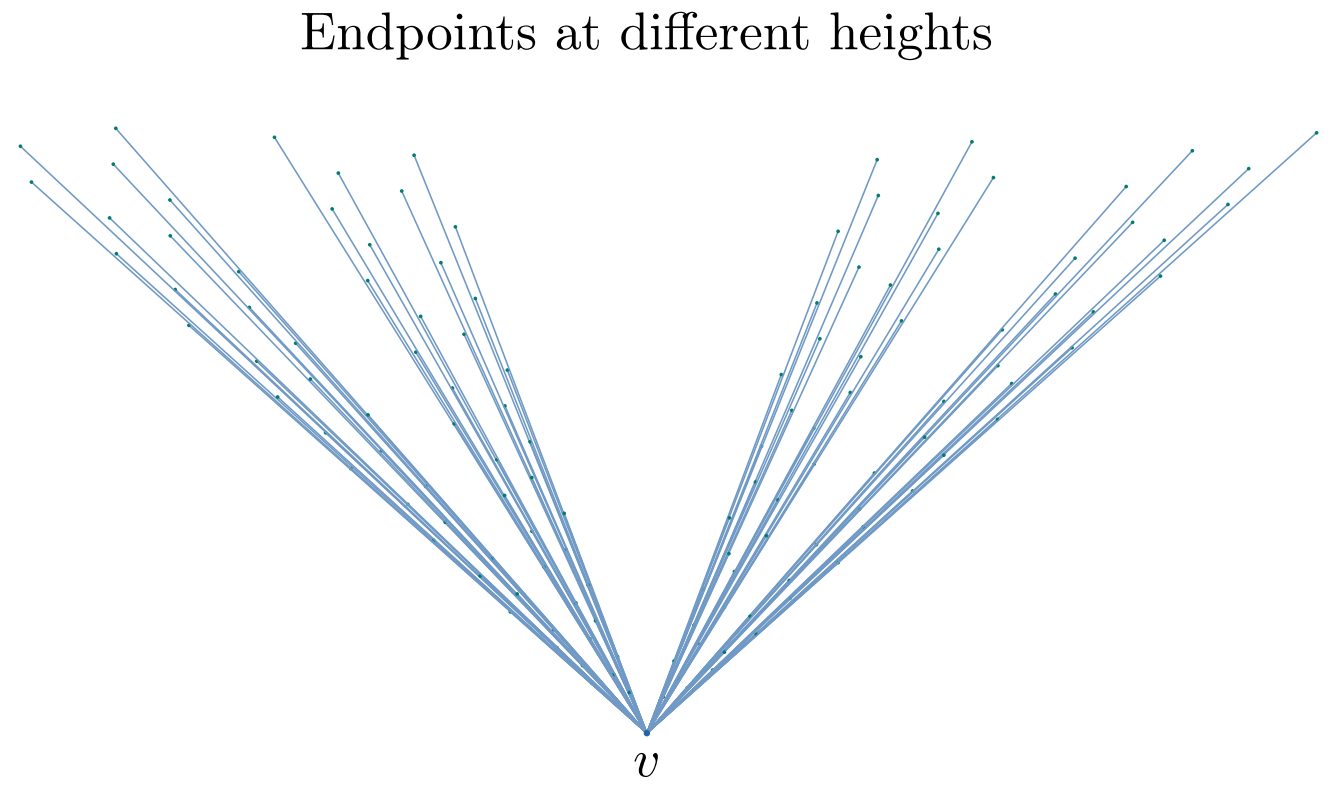  
\end{enumerate}
\end{definition}
\begin{observation}
A fan is smooth if and only if it is a smooth dendroid.	
\end{observation}

For a fan that is an arc, we may choose an interior point as its
top; it then has two legs. For a singleton, its only point is its
top and it has no nondegenerate leg. These conventions only concern
the degenerate and arc cases, and do not affect the assertion about
uncountably many legs.

Note that a Lelek fan was constructed by A.~Lelek in \cite{lelek}.  An interesting  property of the Lelek fan $L$ is the fact that the set of its end points is a dense one-dimensional set in $L$. It is also unique, i.e., any two Lelek fans are homeomorphic, for the proofs see \cite{oversteegen} and \cite{charatonik}.
\begin{definition}
	We  use $\mathcal{SF}$ to denote the class of all smooth fans.
\end{definition}

Recently, the authors proved the following theorem.
\begin{theorem}\label{LelekToCantor}
	There is a continuous surjection from the Lelek fan to the Cantor fan.
\end{theorem}
\begin{proof}
	The authors proved in \cite[Theorem 1.2]{BE} that there is a surjective function $f$ on a Cantor fan $X$ such that $\varprojlim(X,f)$ is a Lelek fan. Then any projection from $\varprojlim(X,f)$ onto $X$ is a continuous surjection from Lelek fan to Cantor fan. 
\end{proof}
\begin{observation}
	Theorem \ref{LelekToCantor} also follows from the earlier work of D. Barto\v sov\' a and A. Kwiatkowska \cite[Theorem 2.12]{bartosova}.  	
\end{observation}
\begin{theorem}\label{BK}
	Let $L$ be a Lelek fan. For each $X\in \mathcal{SF}$, there is a continuous surjection from $L$ to $X$.
\end{theorem}
\begin{proof}
	The theorem follows from Theorems \ref{CE} and \ref{LelekToCantor}.
\end{proof}
\begin{observation}\label{cica}
	Note that Lelek fan is a smooth dendroid. Therefore, it follows from Theorem \ref{CE} that there is a continuous surjection from the Cantor fan onto the Lelek fan. 
\end{observation}
\begin{corollary}\label{lelekova}
	Let $L$ be a Lelek fan and let $C$ be the Cantor fan. Then $L$ and $C$ are both universal domains in $\mathcal{SF}$. 
\end{corollary}
\begin{proof}
	The corollary follows directly from Theorem \ref{BK} and Observation \ref{cica}.
\end{proof}

Theorem \ref{SF} is the main result of this section. 
Several implications in it follow from known
results. The implication \ref{uv3}$\Rightarrow$\ref{uv5} follows
from the retraction theorem of Charatonik and Charatonik
\cite[Theorem~6, p.~165]{CharatonikCharatonik1989}; see also
\cite[Theorem~2.9]{BanicErcegGreenwoodKennedy2025}.
The equivalence of \ref{uv3} and \ref{uv4} follows from the
known mutual embeddability of the Cantor and Lelek fans;
see \cite{BanicErcegKennedy2025}.
The equivalences among \ref{uv1}, \ref{uv5}, \ref{uv6}, and
\ref{uv7} follow from Corollary~\ref{lelekova} and
Theorem~\ref{zelosplosno}. We include explicit arguments
connecting these conditions with the number of legs of $X$.
\begin{theorem}\label{SF}
	Let $L$ be a Lelek fan and let $C$ be a Cantor fan. For each smooth fan $X$, the following statements are equivalent.
	\begin{enumerate}
		\item\label{uv1} $X$ is a reciprocating domain in $\mathcal{SF}$.
		\item\label{uv5} There is a continuous surjection from $X$ to $C$.
		\item\label{uv6} There is a continuous surjection from $X$ to $L$.
		\item \label{uv3} $X$ contains a Cantor fan.
		\item\label{uv4} $X$ contains a Lelek fan.
		\item\label{uv2} $X$ has uncountably many legs.
		\item\label{uv7} $X$ is a universal domain in $\mathcal{SF}$.
	\end{enumerate}
\end{theorem}

\begin{proof}
By Corollary~\ref{lelekova} and Theorem~\ref{zelosplosno},
statements \ref{uv1}, \ref{uv5}, \ref{uv6}, and \ref{uv7} are
equivalent. We prove the remaining implications explicitly.

Suppose that $X$ has uncountably many legs. The endpoint set of a
smooth fan is Polish; see \cite{lelek,HernandezGutierrezHoehn2023}.
By \cite[Corollary 6.5]{Kechris1995}, choose a Cantor set
$M\subseteq E(X)$. Use the standard embedding of a smooth fan
into the Cantor fan, with common top $v$; see
\cite[Section 2.1 and Lemma 2.1]{HernandezGutierrezHoehn2023}.
Write this ambient fan as
\[
H=(K\times[0,1])/(K\times\{0\}),
\]
where $K$ is a Cantor set, and denote its quotient map by $q$.
Away from $v$, both the direction and the height coordinates are
continuous. Since $M$ is compact, avoids $v$, and contains at most
one endpoint on each leg, its set $N\subseteq K$ of directions is
a Cantor set and its height function is positive and continuous
on $N$. Choose $t_0>0$ no larger than its minimum. Since $X$
contains the arc from $v$ to every point of $M$,
\[
F_0=q(N\times[0,t_0])\subseteq X.
\]
The continuum $F_0$ is a Cantor fan. This proves
\ref{uv2}$\Rightarrow$\ref{uv3}.

Suppose next that $X$ contains a Cantor fan $F$. Its top is the
unique ramification point of $X$. The compact set of endpoints
of $F$ gives, by the same argument, a Cantor set $N$ of directions
and a number $t_0>0$ such that $F_0=q(N\times[0,t_0])\subseteq X$.
A non-empty closed subset of the Cantor set is a retract of the
Cantor set; choose a retraction $\rho:K\to N$ (see
\cite[Theorem 7.3]{Kechris1995}). The continuous map
\[
(c,t)\longmapsto q(\rho(c),\min\{t,t_0\})
\]
on $K\times[0,1]$ is constant on $K\times\{0\}$. It induces
a continuous map $R:H\to F_0$ fixing $F_0$ pointwise.
Thus $R|_X$ is a retraction onto the Cantor fan $F_0$.
Composing with a homeomorphism $F_0\to C$ proves
\ref{uv3}$\Rightarrow$\ref{uv5}.

To prove \ref{uv5}$\Rightarrow$\ref{uv2}, suppose instead that
$X$ has at most countably many legs. Then $X$ is a countable
union of arcs and possibly a singleton. For any continuous path
$p:[0,1]\to C$, the set $p^{-1}(C\setminus\{v_C\})$ is an
open subset of $[0,1]$ and has at most countably many components.
The image of each component lies in a single leg of $C$, since
the direction projection into the Cantor set is continuous there.
Thus a path in $C$ meets at most countably many legs away from
$v_C$. It follows that a continuous image of $X$ in $C$ also
meets at most countably many legs away from $v_C$, and cannot
be all of $C$. This proves the desired implication.

Finally, every smooth fan embeds in a Cantor fan, and a Lelek fan
has uncountably many endpoints and hence contains a Cantor fan
by the argument above (also, see \cite{BanicErcegKennedy2025}). Therefore a Cantor fan contains a Lelek
fan and a Lelek fan contains a Cantor fan. This proves the
equivalence of \ref{uv3} and \ref{uv4}, completing the proof.
\end{proof}

\subsection{Indecomposable continua}\label{sec:indecomposable}

We prove that the class of indecomposable continua is directed
upwards, and hence has no reciprocating domain. The common-domain
construction is given by successive doubles of irreducible continua.

\begin{definition}
We use $\mathcal{IC}$ to denote the class of indecomposable continua.
A continuum $X$ is \emph{irreducible between} $a,b\in X$ if no proper
subcontinuum of $X$ contains both $a$ and $b$. It is irreducible
if such points exist.
\end{definition}
The following result is a consequence of a stronger theorem
of Bellamy \cite[Corollary~3]{Bellamy1971}, which states that
there is no continuum admitting continuous surjections onto
all indecomposable metric continua. In particular, there is
no universal element in the class of indecomposable continua
with respect to continuous surjections. We include an
alternative proof in the framework developed here.
\begin{theorem}\label{fujf}
There is no universal domain in $\mathcal{IC}$.
\end{theorem}
\begin{proof}
Ingram constructed a family of indecomposable tree-like continua
such that no continuum maps onto every member of the family;
see \cite{ingram}. In particular, no member of $\mathcal{IC}$
can be a universal domain in $\mathcal{IC}$.
\end{proof}

\begin{lemma}\label{cmok2}
Let $Y,Z$ be subcontinua of a continuum $W$ such that
$W=Y\cup Z$ and $Y\cap Z=\{b\}$. If $Y$ is irreducible
between $a$ and $b$, then every subcontinuum $K$ of $W$
containing $a,b$ contains $Y$.
\end{lemma}
\begin{proof}
The map $R:W\to Y$ which fixes $Y$ and sends $Z$ to $b$ is
continuous. Since $b\in K$, we have $R(K)=K\cap Y$.
This is a subcontinuum of $Y$ containing $a,b$, so it is $Y$.
\end{proof}

\begin{lemma}\label{cmok1}
Let $Y$ be irreducible between $a,b$, let $Z$ be irreducible
between $b,c$, and suppose $Y\cap Z=\{b\}$. Then $Y\cup Z$
is irreducible between $a,c$.
\end{lemma}
\begin{proof}
Let $K\subseteq Y\cup Z$ be a subcontinuum containing $a,c$.
If $b\notin K$, the sets $K\cap Y$ and $K\cap Z$ would
separate $K$, unless one of $a,c$ already equals $b$, in which
case $b\in K$ directly. Thus $b\in K$.
Applying Lemma~\ref{cmok2} to $Y$ and $Z$, we obtain
$Y\subseteq K$ and $Z\subseteq K$. Hence $K=Y\cup Z$.
\end{proof}

The following lemma gives a known criterion for the
indecomposability of an inverse limit in terms of its bonding
maps; see \cite[Theorem~2.7]{nadler} and
\cite[Theorem~2.3]{Eberhart2004}. We include a direct proof
for completeness.
\begin{lemma}\label{Licka}
Let $(Z_n,p_n)$ be an inverse sequence of continua with surjective
bonding maps. Suppose that for every positive integer $n$ and
all subcontinua $H,K$ of $Z_{n+1}$,
\[
Z_{n+1}=H\cup K\quad\Longrightarrow\quad
p_n(H)=Z_n\ \text{ or }\ p_n(K)=Z_n.
\]
Then $\invlim(Z_n,p_n)$ is indecomposable.
\end{lemma}
\begin{proof}
Write $Z=\invlim(Z_n,p_n)$ and let $\pi_n$ be its coordinate
projections. A proper subcontinuum $A$ of $Z$ has a proper
coordinate image. Indeed, if $\pi_n(A)=Z_n$ for every $n$,
then for each $z=(z_n)\in Z$ the nested compact sets
\[
\{a\in A:\pi_n(a)=z_n\},\qquad n\in\mathbb N,
\]
are non-empty and have non-empty intersection. Their intersection
can contain only $z$, so $A=Z$, a contradiction.

Suppose $Z=A\cup B$ for proper subcontinua $A,B$. Choose $n$
large enough that both $\pi_n(A)$ and $\pi_n(B)$ are proper in
$Z_n$. Such a common $n$ exists, since a proper image at one
coordinate forces proper images at every later coordinate.
Now $Z_{n+1}=\pi_{n+1}(A)\cup\pi_{n+1}(B)$, so the hypothesis
implies $\pi_n(A)=Z_n$ or $\pi_n(B)=Z_n$, a contradiction.
\end{proof}
The following theorem is a special case of a stronger result
of Bellamy \cite[Corollary~4]{Bellamy1971}, which states that
every metric continuum is a retract of an indecomposable
metric continuum. Indeed, for any two continua $A$ and $B$,
their wedge $A\vee B$ is a continuum admitting retractions
onto both $A$ and $B$. Bellamy's result therefore yields
an indecomposable continuum containing copies of $A$ and
$B$ as retracts. In particular, it admits continuous
surjections onto both continua, without any assumption
of irreducibility. We include a direct proof using inverse
limits, which provides an explicit construction in the
irreducible case.

\begin{theorem}\label{ired}
For any two irreducible continua $A,B$, there is an indecomposable
continuum $X$ containing homeomorphic copies of $A$ and $B$ and
admitting continuous surjections onto both $A$ and $B$.
\end{theorem}
\begin{proof}
If both continua are singletons, take $X$ to be a singleton.
Otherwise, let $A$ be irreducible between $a,c$, and let $B$
be irreducible between $b,d$. Form $Z_1$ by identifying $c$ with
$b$ in the disjoint union of $A$ and $B$. Let $a_1,d_1$ be the
images of $a,d$. By Lemma~\ref{cmok1}, $Z_1$ is a nondegenerate
continuum irreducible between $a_1,d_1$.

Suppose that $Z_n$ has been constructed and is irreducible
between distinct points $a_n,d_n$. Take two copies $Z_n^0,Z_n^1$
of $Z_n$ and identify their copies of $d_n$. Let $Z_{n+1}$ be
the resulting wedge, let $c_{n+1}$ be the identified point, and
let $a_{n+1},d_{n+1}$ be the copies of $a_n$ in $Z_n^0,Z_n^1$,
respectively. By Lemma~\ref{cmok1}, $Z_{n+1}$ is irreducible
between $a_{n+1},d_{n+1}$. Define
$p_n:Z_{n+1}\to Z_n$ to equal the canonical homeomorphism onto
$Z_n$ on each copy. These definitions agree at $c_{n+1}$, so
$p_n$ is a continuous surjection. Let
$i_n:Z_n\to Z_{n+1}$ be the inclusion into the first copy.
Then $p_n\circ i_n=\id_{Z_n}$.

We check the hypothesis of Lemma~\ref{Licka}. If
$Z_{n+1}=H\cup K$, one of $H,K$ contains at least two of
$a_{n+1},c_{n+1},d_{n+1}$. Suppose it is $H$. If $H$ contains
$a_{n+1},c_{n+1}$, Lemma~\ref{cmok2} implies $Z_n^0\subseteq H$.
If it contains $c_{n+1},d_{n+1}$, the same lemma gives
$Z_n^1\subseteq H$. If it contains both $a_{n+1},d_{n+1}$,
irreducibility implies $H=Z_{n+1}$. In every case,
$p_n(H)=Z_n$. The case with $K$ in place of $H$ is identical.
Thus $X=\invlim(Z_n,p_n)$ is indecomposable.

The map
\[
z\longmapsto
\bigl(z,i_1(z),i_2(i_1(z)),\ldots\bigr)
\]
is a continuous injection $Z_1\to X$, and hence an embedding.
Therefore $X$ contains copies of $A$ and $B$. Finally, the first
projection $X\to Z_1$ is surjective. Composing it with the
retractions of $Z_1$ onto the copies of $A$ and $B$ gives
continuous surjections from $X$ onto $A$ and $B$.
\end{proof}

\begin{lemma}\label{ir}
Every indecomposable continuum is irreducible.
\end{lemma}
\begin{proof}
For a nondegenerate indecomposable continuum, take two points
from distinct composants; no proper subcontinuum contains both.
A singleton is irreducible between its only point and itself.
\end{proof}
The following theorem is also an immediate consequence of
Bellamy's result \cite[Corollary~4]{Bellamy1971}.
We give a short alternative proof using Theorem~\ref{ired}.
\begin{theorem}\label{d2va}
For any two indecomposable continua $A,B$, there is an
indecomposable continuum containing copies of both and mapping
continuously onto both.
\end{theorem}
\begin{proof}
Apply Lemma~\ref{ir} and Theorem~\ref{ired}.
\end{proof}

\begin{corollary}\label{fuj}
The class $\mathcal{IC}$ is directed upwards.
\end{corollary}
\begin{proof}
This follows from Theorem~\ref{d2va}.
\end{proof}

\begin{theorem}
There are no reciprocating domains in $\mathcal{IC}$.
\end{theorem}
\begin{proof}
By Corollary~\ref{fuj} and Theorem~\ref{cilka}, a reciprocating
domain in $\mathcal{IC}$ would be a universal domain in
$\mathcal{IC}$, contrary to Theorem~\ref{fujf}.
\end{proof}

\section*{Acknowledgments}

The authors acknowledge the use of artificial intelligence tools in developing the construction and proofs in Subsubsections~\ref{sec:pa-construction-start}--\ref{pa-sec:main}, which establish the existence of a chainable continuum mapping onto a pseudoarc and containing no pseudoarc.
Work on this construction began at the 40th Summer Conference on Topology and its Applications, held in Split, Croatia, on July 13--17, 2026, with the question of whether every chainable continuum that maps continuously onto a pseudoarc must contain a pseudoarc.
The authors provided the AI tools with specific mathematical ideas and directions for constructing a counterexample. In particular, they proposed modifying an inverse-limit representation of the pseudoarc by replacing certain parts of the coordinate intervals with nondegenerate intervals and appropriately defined bonding maps, with the aim of introducing a dense collection of arcs while preserving chainability and the existence of a continuous surjection onto the pseudoarc. Guided by these ideas, the AI tools developed the technical construction of the arcs using methods and techniques already known in continuum theory.
The construction and its proofs were subsequently developed through successive revisions, including identifying gaps, checking degenerate cases, and refining the lifting argument. The AI tools also assisted with the exposition and illustrations of the construction.
To assess the role of the mathematical guidance provided by the authors, the same problem was independently presented to an AI tool through a separate user account, without any additional suggestions or instructions concerning a possible construction. In that instance, the AI tool was unable to solve the problem. This comparison illustrates the importance of the authors' mathematical ideas and specific guidance in obtaining the construction presented in this paper.
The authors take full responsibility for the mathematical correctness of the results and the final content of the paper.

\noindent I. Bani\v c\\
              (1) Faculty of Natural Sciences and Mathematics, University of Maribor, Koro\v{s}ka 160, SI-2000 Maribor,
   Slovenia; \\(2) Institute of Mathematics, Physics and Mechanics, Jadranska 19, SI-1000 Ljubljana, 
   Slovenia; \\(3) Andrej Maru\v si\v c Institute, University of Primorska, Muzejski trg 2, SI-6000 Koper,
   Slovenia\\
             {iztok.banic@um.si}          
     
				\-
				
		\noindent G.  Erceg\\
             Faculty of Science, University of Split, Rudera Bo\v skovi\' ca 33, Split,  Croatia\\
    
{{goran.erceg@pmfst.hr}       }             

                 \-
	
  \noindent J.  Kennedy\\
             Lamar University, 200 Lucas Building, P.O. Box 10047, Beaumont, TX 77710 USA\\
   
{{kennedy9905@gmail.com}       } 

              \-
	
  \noindent C.  Mouron\\
             Rhodes College,  2000 North Parkway, Memphis, Tennessee 38112  USA\\\
     
{{mouronc@rhodes.edu}       } 

              \-
	
  \noindent V.  Nall\\
            Department of Mathematics,  University of Richmond, Richmond, Virginia 23173 USA\\
  
{{vnall@richmond.edu}       }    
                 
     \end{document}